\documentclass[10pt]{article}
\usepackage{graphicx}
\usepackage{epsf}
\usepackage{epsfig}

\usepackage{amsmath}
\usepackage{amsfonts}
\usepackage{amssymb}

\newcommand\ve{\varepsilon}
\newcommand\ra{\rightarrow}

\newcommand\cT{\mathcal T}
\newcommand\bs{\backslash}

\def\lsim{\hbox{\kern -.2em\raisebox{-1ex}{$~\stackrel{\textstyle<}{\sim}~$}}\kern -.2em}
\def\gsim{\hbox{\kern -.2em\raisebox{-1ex}{$~\stackrel{\textstyle>}{\sim}~$}}\kern -.2em}

\def\sq{\hfill $\diamond$\\}
\def\be{\begin{equation}}
\def\ee{\end{equation}}

\newcommand\iref[1]{(\ref{#1})}
\def\cL{\mathcal L}
\newcommand\R{\mathbb{R}}

\newcommand\proof{\paragraph{Proof: }}

\DeclareMathOperator\Id{Id}
\DeclareMathOperator\disc{disc}

\newtheorem{definition}{Definition}[section]
\newtheorem{theorem}{Theorem}[section]
\newtheorem{prop}{Proposition}[section]
\newtheorem{lemma}[prop]{Lemma}
\newtheorem{corol}{Corollary}[section]

\newcommand\mT{\mathbb T}

\newcommand\mi{\mathbf i}

\newcommand\cE{{\mathcal E}}

\def\I{{\rm \hbox{I\kern-.2em\hbox{I}}}}
\def\P{{\rm \hbox{I\kern-.2em\hbox{P}}}}
\def\H{{\rm \hbox{I\kern-.2em\hbox{H}}}}
\def\R{{\rm \hbox{I\kern-.2em\hbox{R}}}}
\def\N{{\rm \hbox{I\kern-.2em\hbox{N}}}}
\def\Z{{\rm {{\rm Z}\kern-.28em{\rm Z}}}}
\def\C{{\rm \hbox{C\kern -.5em {\raise .32ex \hbox{$\scriptscriptstyle
|$}}\kern-.22em{\raise .6ex \hbox{$\scriptscriptstyle |$}}\kern .4em}}}

\def\cT{{\cal T}}

\def\cR{{\cal R}}

\def\cL{{\cal L}}

\def\cP{{\cal P}}

\def\cO{{\cal O}}
\def\Chi{\raise .3ex
\hbox{\large $\chi$}} 
\def\lsima{\hbox{\kern -.6em\raisebox{-1ex}{$~\stackrel{\textstyle<}{\sim}~$}}\kern -.4em}
\def\lsim{\hbox{\kern -.2em\raisebox{-1ex}{$~\stackrel{\textstyle<}{\sim}~$}}\kern -.2em}
\def\({\Bigl (}
\def\){\Bigr )}
\def\({\Bigl (}
\def\){\Bigr )}
\DeclareMathOperator{\diam}{diam}
\DeclareMathOperator{\diag}{diag}
\DeclareMathOperator{\bary}{bary}
\DeclareMathOperator{\cyc}{\mathbb S}

\def\<{\langle}
\def\>{\rangle}
\def\t{\tilde}

\def\e{\varepsilon}

\def\nl{\newline}

\def\b{\mathbf}

\newtheorem{remark}[theorem]{Remark}

\newcommand\seqT{(\cT_N)_{N\geq N_0}}

\newcommand\mhalf{{s_m}}
\newcommand\Teq{{T_\text{eq}}}

\newcommand\trans{\mathrm T}

\newcommand\cTr{{\cT_s^\text{reg}}}
\newcommand\cTb{{\cT_s^\text{bd}}}

\newcommand\Kpol{\b K}

\newcommand\NF{N\hspace{-0.8mm}F}
\newcommand\sep{\; : \;}

\newcommand\Keq{K_{\rm eq}}

\title{Optimal meshes for finite elements of arbitrary order}
\author{Jean-Marie Mirebeau}

\begin{document}

\maketitle

\begin{abstract}
Given a function $f$ defined on a bounded domain $\Omega\subset \R^2$ 
and a number $N>0$, we study the
properties of the triangulation $\cT_N$ that minimizes
the distance between $f$ and its interpolation
on the associated finite element space,
over all triangulations of at most $N$ elements.
The error is studied in the norm $X=L^p$ for $1\leq p\leq \infty$ and we consider Lagrange finite
elements of arbitrary polynomial degree $m-1$.
We establish sharp asymptotic error estimates as $N\to +\infty$ when the optimal
anisotropic triangulation is used, recovering the results on piecewise linear interpolation \cite{BBLS,B,CSX}, an improving the results on higher degree interpolation \cite{C3,C1,C2}.
These estimates involve invariant polynomials applied to the $m$-th order
derivatives of $f$. In addition, our analysis also provides with practical strategies
for designing meshes such that the interpolation error satisfies the optimal estimate 
up to a fixed multiplicative constant. We partially extend our results to higher dimensions
for finite elements on simplicial partitions of a domain $\Omega\subset \R^d$.
\vspace{-0.1cm}
\vspace{-0.5cm}
\end{abstract}

\section{Introduction.}

\subsection{Optimal mesh adaptation}

In finite element approximation, a usual distinction is between {\it uniform} and
{\it adaptive} methods. In the latter, the elements defining the mesh 
may vary strongly in size and shape for a better adaptation
to the local features of the approximated function $f$. This naturally raises
the objective of characterizing and constructing an {\it optimal mesh}
for a given function $f$. 

Note that depending on the context, the function $f$ may be fully
known to us, either through an explicit formula or a discrete sampling, 
or observed through noisy measurements,
or implicitly defined as the solution of a given
partial differential equation.

In this paper, we assume that $f$ is a function 
defined on a polygonal bounded domain $\Omega\subset \R^2$.
For a given conforming triangulation $\cT$ of $\Omega$,
and an arbitrary but fixed 
integer $m>1$, we denote by $I_{m,\cT}$ 
the standard interpolation operator on the 
Lagrange finite elements of degree $m-1$
space associated to $\cT$.
Given a norm $X$ of interest and a number
$N>0$, the objective of finding the optimal mesh for $f$ can be
formulated as solving the optimization problem
$$
\min_{\#(\cT)\leq N} \|f-I_{m,\cT}f\|_X,
$$
where the minimum is taken over all conforming triangulations of 
cardinality $N$. We denote by $\cT_N$ the minimizer of the
above problem.

Our first objective is to establish
sharp asymptotic error estimates that precisely describe the
behavior of $\|f-I_{m,\cT}f\|_X$ as $N\to +\infty$. Estimates of that type
were obtained in \cite{BBLS,B,CSX} in the particular case of linear
finite elements ($m-1=1$) and with the error measured in $X=L^p$. They have the form
\be
\limsup_{N\to +\infty} \(N\min_{\#(\cT)\leq N}\|f-I_{m,\cT}f\|_{L^p}\)  \leq C\| \sqrt{|\det(d^2f)|}\|_{L^\tau},\;\; \frac 1 \tau=\frac 1 p+1,
\label{optiaffine}
\ee
which reveals that the convergence rate is governed by the quantity $\sqrt{|\det(d^2f)|}$, which depends
nonlinearly the Hessian $d^2f$. This is heavily tied to the fact that we allow
triangles with possibly highly anisotropic shape.
In the present work, the polynomial degree $m-1$ is arbitrary
and the quantities governing the convergence rate will therefore depend nonlinearly on the 
$m$-th order derivative $d^mf$.

Our second objective is to propose simple and practical ways
of designing meshes which behave similar to the optimal one,
in the sense that they satisfy the sharp error estimate up
to a fixed multiplicative constant.

\subsection{Main results and layout}

We denote by $\H_m$ the space of homogeneous polynomials of degree $m$
$$
\H_m:={\rm Span}\{ x^ky^l\sep k+l=m\}.
$$
For any triangle $T$, we denote by $I_{m,T}$ the local interpolation
operator acting from $C^0(T)$ onto $\P_{m-1}$ the space of 
polynomials of total degree $m-1$. The image of $v\in C^0(T)$ by this operator is defined by the
conditions
$$
I_{m,T}v(\gamma)=v(\gamma),
$$
for all points $\gamma\in T$ with barycentric coordinates in
the set $\{0, \frac 1 {m-1},\frac 2 {m-1},\cdots,1\}$. We denote by
$$
e_{m,T}(v)_p:=\|v-I_{m,T}v\|_{L^p(T)}
$$
the interpolation error measured in the norm $L^p(T)$. We also denote
by
$$
e_{m,\cT}(v)_p:=\|v-I_{m,\cT}v\|_{L^p}=\left(\sum_{T\in\cT}e_{m,T}(v)_p^p\right)^{\frac 1 p},
$$
the global interpolation error for a given triangulation $\cT$, with the standard modification if $p=\infty$.

A key ingredient in this paper is a function 
defined by a {\it shape optimization problem}:
for any fixed $1\leq p\leq \infty$ and for any $\pi\in \H_m$, we define
\be
K_{m,p}(\pi):=\inf_{|T|=1}e_{m,T}(\pi)_p.
\label{shapefunction}
\ee
Here, the infimum is taken over all triangles of area $|T|=1$. 
Note that from the homogeneity of $\pi$, we find that
\be
\inf_{|T|=A}e_{m,T}(\pi)_p=K_{m,p}(\pi)A^{\frac m 2 + \frac 1 p}.
\label{shapeA}
\ee
This optimization problem thus gives the shape of the triangles 
of a given area which is at best adapted
to the polynomial $\pi$ in the sense of minimizing the interpolation error
measured in $L^p$.  
We refer to $K_{m,p}$ as the {\it shape function}.
We discuss in \S 2 the main properties of this function. 
\nl
\nl
Our asymptotic error estimate for the optimal triangulation is given by
the following theorem.
\begin{theorem}
\label{maintheorem}
For any polygonal domain $\Omega\subset \R^2$, and any function 
$f\in C^m(\Omega)$, there exists a sequence 
of triangulations $\seqT$, with $\#(\cT_N)=N$ such that
\be
\label{Cf}
\limsup_{N\ra \infty} N^{\frac m 2} e_{m,\cT_N}(f)_p\leq 
\left\|K_{m,p}\left(\frac{d^m f}{m!}\right)\right\|_{L^q(\Omega)},\ \frac 1 q := \frac m 2 +\frac 1 p 
\ee
\end{theorem}
An important feature of this estimate is the ``$\limsup$'' asymptotical operator. Recall that the upper limit of a sequence $(u_N)_{N\geq N_0}$ is defined by 
$$
\limsup_N u_N := \lim_{N\to \infty} \sup_{n\geq N} u_n,
$$
and is in general strictly smaller than the supremum $\sup_{N\geq N_0} u_N$. It is still an open question to find an appropriate upper estimate for $\sup_N N^{m/2} e_{m,\cT_N}(f)_p$ when optimally adapted anisotropic triangulations are used.

In the estimate \iref{Cf}, the $m$-th derivative $d^mf$ is identified to 
an homogeneous polynomial in $\H_m$:
$$
\frac{d^mf}{m!}\sim \sum_{k+l=m}\frac{\partial^m f}{\partial^kx\partial^l y}\frac{x^k}{k!}\frac{y^l}{l!}.
$$
In order to illustrate the sharpness of \iref{Cf}, we introduce
a slight restriction on sequences of triangulations, following 
an idea in \cite{BBLS}: a sequence $\seqT$ of triangulations, such that $\#(\cT_N) = N$, is said to be \emph{admissible} if
\be
\label{admissibilitycond}
\sup_{T\in \cT_N} \diam(T) \leq C_AN^{-1/2},
\ee
for some $C_A>0$ independent of $N$. The following theorem shows that the estimate
\iref{Cf} cannot be improved when we restrict our attention to admissible sequences.
It also shows that this class is reasonably large in the sense that 
\iref{Cf} is ensured to hold up to small perturbation.

\begin{theorem}
\label{optitheorem}
Let $\Omega\subset \R^2$ be a compact polygonal domain, and $f\in C^m(\Omega)$. 
Denote $ \frac 1 q := \frac m 2 +\frac 1 p$.
For all \emph{admissible} sequences of triangulations $\seqT$, one has
$$
\liminf_{N\ra \infty} N^{\frac m 2} e_{m,\cT_N}(f)_p \geq \left\|K_{m,p}\left(\frac{d^m f}{m!}\right)\right\|_{L^q(\Omega)}.
$$
For all $\ve>0$, there exists an \emph{admissible} sequence of triangulations $(\cT_N^\ve)_{N\geq N_0}$, such that
$$
\limsup_{N\ra \infty} N^{\frac m 2} e_{m,\cT_N^\ve}(f)_p \leq \left\|K_{m,p}\left(\frac{d^m f}{m!}\right)\right\|_{L^q(\Omega)}+\ve.
$$
\end{theorem}
Note that the sequences $(\cT_N^\ve)_{N\geq N_0}$ satisfy the admissibility condition \iref{admissibilitycond} with a constant $C_A(\ve)$ which may explode as $\ve\to 0$.
The proofs of both theorems are given in \S 3. These proofs reveal that the construction of the 
optimal triangulation obeys two principles:
(i) the triangulation should {\it equidistribute} the local approximation error $e_{m,T}(f)_p$ between 
each triangle and (ii) the aspect ratio of a triangle $T$ should be 
{\it isotropic} with respect to a distorted
metric induced by the local value of $d^mf$ on $T$
(and therefore anisotropic in the sense of the euclidean metric). 
Roughly speaking, the quantity $\|K_{m,p}\left(\frac{d^m f}{m!}\right)\|_{L^q(T)}$
controls the local interpolation $L^p$-error estimate on a triangle $T$ once this triangle  
is optimized with respect to the local properties of $f$.
This type of estimate differs from those obtained in \cite{Ap} which hold for any $T$, optimized
or not, and involve the partial derivatives of $f$ in a local coordinate system which is 
adapted to the shape of $T$.

The proof of the upper estimates in Theorem \ref{optitheorem} involves the construction
of an optimal mesh based on a patching strategy similar to \cite{B}. 
However, inspection of the proof reveals that
this construction becomes effective only when 
the number of triangles $N$ becomes very large. Therefore it 
may not be useful in practical applications.

A more practical approach consists in deriving the above mentioned distorted metric from 
the exact or approximate data of $d^mf$, using the following procedure. 
To any $\pi\in \H_m$, we associate a 
symmetric positive definite matrix $h_\pi\in S_2^+$. If $z\in \Omega$ and $d^m f(z)$ is close to $\pi$, then the triangle $T$ containing $z$ should be isotropic in the metric $h_\pi$.
The global metric is given
at each point $z$ by
$$
h(z)=s(\pi_z) h_{\pi_z},\;\; \pi_z=d^mf(z),
$$
where $s(\pi_z)$ is a scalar factor which depends on the desired accuracy
of the finite element approximation. 
Once this metric has been properly identified, fast algorithms such as in \cite{Inria, Bamg, Peyre}
can be used to design a near-optimal mesh based on it. Recently in \cite{Shew,Bois}, several algorithms have been rigorously proved to terminate and produce good quality meshes.
Computing the map
\be
\pi\in \H_m \mapsto h_\pi\in S_2^+,
\label{mappi}
\ee
is therefore of key use in applications. This problem is well understood
in the case of linear elements ($m=2$): the matrix $h_\pi$ is then defined as the
absolute value (in the sense of symmetric matrices) of the matrix associated to
the quadratic form $\pi$. In contrast, the exact form of this map in the case $m\geq 3$
is not well understood. 

In this paper, we propose
algebraic strategies for computing the map \iref{mappi} for $m=3$ which corresponds to
quadratic elements. These strategies
have been implemented in an open-source 
Mathematica code \cite{sitejm}. 
In a similar manner, we address the algebraic computation 
of the shape function $K_{m,p}(\pi)$
from the coefficients of $\pi\in \H_m$, when $m\geq 3$. 
All these questions are addressed in \S 4, 5 and 6.

In \S 4, we discuss the particular case of linear ($m=2$)
and quadratic ($m=3$) elements. In this case, it is possible
to obtain explicit formulas for $K_{m,p}(\pi)$ from the coefficients
of $\pi$. In the case $m=2$, this formula is of the form
$$
K_{2,p}(ax^2+2bxy+cy^2)=\sigma\sqrt{|b^2-ac|},
$$
where the constant $\sigma$ only depends on $p$ and the sign of $b^2-ac$,
and we therefore recover the known estimate \iref{optiaffine} from Theorem \ref{maintheorem}. The formula for $m=3$ involves the
discriminant of the third degree polynomial $d^3f$. Our analysis
also leads to an algebraic computation of the map \iref{mappi}.
We want to mention that a different strategy for the
the construction of the distorted metric and
the derivation of error estimate for finite element
of arbitrary order was proposed in
\cite{C3}. In this approach, the distorted
metric is obtained at a point $z\in\Omega$ by finding the largest
ellipse contained in a level set of the polynomial $d^mf_z$.
This optimization problem has connections with the one
that defines the shape function in \iref{shapefunction}
as we shall explain in \S 2. The approach 
proposed in the present work in the case $m=3$ has the advantage of
avoiding the use of numerical optimization,
the metric being directly derived from the
coefficients of $d^mf$. 

In \S 5, we address the case $m>3$. In this case, explicit formulas for $K_{m,p}(\pi)$
seem out of reach. However we can introduce explicit functions 
$\Kpol_m(\pi)$ which are polynomials in the coefficients of $\pi$, 
and are equivalent to $K_{m,p}(\pi)$, leading therefore to similar
asymptotic error estimates up to multiplicative constants. At the current stage,
we did not obtain a simple solution to the
algebraic computation of the map \iref{mappi} in the case $m>3$.
The derivation of $\Kpol_m$ is based on the 
theory of invariant polynomials due to Hilbert. Let us mention
that this theory was also recently applied in \cite{OST} to image processing tasks
such as affine invariant edge detection and denoising.

We finally discuss in \S 6 the possible extension of our analysis to simplicial
elements in higher dimension. This extension is not straightforward except in the case of
linear elements $m=2$.

\section{The shape function}

In this section, we establish several properties of the 
function $K_{m,p}$ which will be of key use in the sequel.
We assume that $m\geq 2$ is an integer, and $p\in[1,\infty]$.
We equip the finite dimensional vector space $\H_m$ with a norm $\|\cdot\|$
defined as the supremum of the coefficients
\be
\label{normdef}
\text{ If } \pi(x,y)=\sum_{i=0}^m a_i x^i y^{m-i}, \text{ then } \|\pi\| = \max_{0\leq i\leq m} |a_i|.
\ee
Our first result shows that the function $K_{m,p}$
vanishes on a set of polynomials which has a simple algebraic characterization.
\begin{prop}
\label{vanishprop}
We denote by $\mhalf:= \lfloor \frac m 2 \rfloor +1$ the smallest integer strictly larger than $m/2$.
The vanishing set of $K_{m,p}$ is the set of polynomials which have a generalized 
root of multiplicity at least $\mhalf$:
$$
K_{m,p}(\pi)=0 \Leftrightarrow  \pi(x,y) = (\alpha x +\beta y)^\mhalf \tilde \pi, \mbox{ for some }
\alpha,\beta \in \R\mbox{ and } \t\pi \in \H_{m-\mhalf}.
$$
\end{prop}

\proof
We denote by $\Teq$ a fixed equilateral triangle of unit area, centered at $0$.

We first assume that $\pi(x,y)= (\alpha x+\beta y)^\mhalf \tilde \pi$. Then there exists a rotation  
$R\in \cO_2$ and $\hat \pi\in H_{m-\mhalf}$ such that 
$$
\pi\circ R(x,y)= x^\mhalf \hat \pi(x,y)=x^\mhalf \left( \sum_{i=0}^{m-s_m}a_ix^iy^{m-s_m-i} \right),
$$
Therefore denoting by $\phi_\ve$ the linear transform $\phi_\ve(x,y) = R\left(\ve x ,\frac y \ve\right)$ we obtain
$$
\|\pi\circ \phi_\ve\| =  \max_{i=0,\cdots,m-s_m} |a_i| \ve^{2s_m-m+2i}
\leq \ve^{2\mhalf-m} \|\hat \pi\|
\to 0\; \; {\rm as}\;\; \ve\to 0.
$$
Consequently 
$$
e_{m,\phi_\ve(\Teq)} (\pi)_p = e_{m,\Teq} (\pi\circ \phi_\ve)_p\to 0\; \; {\rm as}\;\; \ve\to 0.
$$
Since $|\det \phi_\ve| = 1$, the triangles $\phi_\ve(\Teq)$ have unit area, 
and therefore $K_{m,p}(\pi) = 0$.

Conversely, let $\pi\in \H_m\bs\{0\}$ be such that $K_{m,p}(\pi)=0$. Then there exists a sequence $(T_n)_{n\geq 0}$ of triangles with unit area such that $e_{m,T_n}(\pi)_p\to 0$.
We remark that the interpolation error $e_T(\pi)_p$ of $\pi\in \H_m$ is 
invariant by a translation $\tau_h: z\mapsto z+h$ of the triangle $T$. Indeed 
$\pi-\pi\circ\tau_h\in \P_{m-1}$ so that
\be
\|\pi- I_{m,T}\pi\|_{L^p(\tau_h(T))}=
\|\pi\circ \tau_h - I_{m,T}(\pi\circ\tau_h)\|_{L^p(T)}=
\|\pi- I_{m,T}\pi\|_{L^p(T)}.
\label{transinv}
\ee
Hence we may assume that the barycenter of $T_n$ is $0$, 
and write $T_n = \phi_n(\Teq)$, for some linear transform $\phi_n$ with $\det \phi_n = 1$.
Since $e_{m,\Teq}(\cdot)_p$ is a norm on $\H_m$, it follows 
that $\pi\circ \phi_n\ra 0$.

The linear transform $\phi_n$ has a singular value decomposition
$$
\phi_n = U_n \circ D_n \circ V_n, \text{ where } U_n,V_n\in \cO_2,\text{ and } D_n = \left(
\begin{array}{cc}
\ve_n & 0\\
0 & 1/\ve_n
\end{array}
\right)
,\ 0<\ve_n\leq 1.
$$
Since the orthogonal group $\cO_2$ is compact, there is a uniform constant $C$ such that 
$$
\|\pi \circ V\| \leq C \|\pi\|,\;\; \pi\in \H_m,\; V\in \cO_2.
$$
Therefore 
$$
\|\pi\circ U_n \circ D_n\| = \|\pi\circ U_n \circ D_n \circ V_n \circ V^{-1}_n\|
\leq C\|\pi\circ\phi_n\| \to 0.
$$
Denoting by $a_{i,n}$ the coefficient of $x^i y^{m-i}$ in
$\pi\circ U_n$, we find that $a_{i,n}\ve_n^{2i-m}$ tends to $0$ as $n\to +\infty$.
In the case where $i<s_m$, this implies that $a_{i,n}$ tends to $0$ as $n\to +\infty$ 

Moreover, again by compactness of $\cO_2$,
we may assume, up to a subsequence, that $U_n$ converges to some $U\in \cO_2$. 
Denoting by $a_i$ the coefficient of  $x^i y^{m-i}$ in
$\pi\circ U$, we thus find that $a_i=0$ if $i<s_m$. This 
This implies that $\pi\circ U (x,y) = x^\mhalf \hat \pi(x,y)$ which concludes the proof.
\sq

\begin{remark}
In the simple case $m=2$, we infer from Proposition \ref{vanishprop}
that $K_{2,p}(\pi)=0$ if and only if $\pi$ is of the form $\pi(x,y)=x^2$ up to a rotation, 
and therefore a one-dimensional function. For such a function, the optimal
 triangle $T$ degenerates to a segment in the $y$ direction, i.e. 
optimal triangles of a fixed area tend to be infinitely long in one direction. 
This situation also holds when $m>2$. Indeed, we see in the second part
in the proof of Proposition \ref{vanishprop} that if $\pi$ is a non-trivial polynomial such that
$K_{m,p}(\pi)=0$, then $\e_n$ must tends to $0$ as $n\to +\infty$. This shows that 
$T_n=\phi_n(T)$ tends to be infinitely flat in the direction $Ue_y$ with $e_y=(0,1)$.
However, $K_{m,p}(\pi)=0$ does not any longer mean that
$\pi$ is a polynomial of one variable.
\end{remark}

\noindent
Our next result shows that the function $K_{m,p}$ is homogeneous, 
and obeys an invariance property with respect to linear change of variables.
\begin{prop}
\label{propinvarK}
For all $\pi \in \H_m$, $\lambda\in \R$ and $\phi \in \cL(\R^2)$,
\begin{eqnarray}
\label{homogK}
K_{m,p}(\lambda \pi) &=& |\lambda| K_{m,p}(\pi)\\
\label{invK}
K_{m,p}(\pi\circ \phi) &=& |\det\phi|^{m/2} K_{m,p}(\pi)
\end{eqnarray}
\end{prop}

\proof
The homogeneity property \iref{homogK} is a direct consequence of the definitions of $K_{m,p}$. In order to prove the invariance property \iref{invK} we assume in a first part that $\det \phi\neq 0$ and
we define $\t T:=\frac{\phi(T)}{\sqrt{|\det \phi|}}$
and $\t \pi (z):=\pi(\sqrt{|\det \phi|}z)= |\det \phi|^{m/2} \pi(z)$. 

We now remark that the local interpolant $I_{m,T}$ commutes 
with linear change of variables in the sense that, when $\phi$ is an invertible linear transform,
\be
I_{m,T}(v\circ \phi)=(I_{m,\phi(T)} v) \circ \phi,
\label{intercommut}
\ee
for all continuous function $v$ and triangle $T$. Using this commutation formula we obtain 
\begin{eqnarray*}
e_{m,T}(\pi\circ \phi)_p &=& |\det \phi|^{-1/p} e_{m,\phi(T)}(\pi)_p\\
&=& e_{m,\t T}(\t \pi)_p\\
&=& |\det \phi|^{m/2} e_{m,\t T} (\pi)_p.
\end{eqnarray*} 

Since the map $T\mapsto \t T$ is a bijection of the set of triangles onto itself, leaving the area invariant, we obtain the relation \iref{invK} when $\phi$ is invertible.
When $\det \phi = 0$, the polynomial $\pi \circ \phi$ can be written $(\alpha x+\beta y)^m$
so that $K_{m,p}(P\circ \phi) = 0$ by Proposition \ref{vanishprop}.
\sq

The functions $K_{m,p}$ are not necessarily continuous, but the 
following properties will be sufficient for our purposes.

\begin{prop}
\label{propsemicont}
The function $K_{m,p}$ is upper semi-continuous in general, and continuous if $m=2$ or $m$ is odd.
Moreover the following property holds:
\be
\text{If } \pi_n\ra \pi \text{ and } K_{m,p}(\pi_n) \ra 0 \text{ then } K_{m,p}(\pi) = 0.
\label{contzero}
\ee
\end{prop}

\proof
The upper semi-continuity property comes from the fact that the infimum of a family of upper
semi-continuous functions is an upper semi-continuous function. 
We apply this fact to the functions $\pi \mapsto e_{m,T}(\pi)_p$ indexed by triangles which are obviously continuous.

For any polynomial $\pi \in \H_2$, $\pi = a x^2+2 b xy+ c y^2$, we define $\det \pi = ac-b^2$. It will be shown in \S4 that $K_{2,p}(\pi)= \sigma_p\sqrt{|\det \pi|}$, where $\sigma_p$ only depends on the sign of $\det \pi$. This clearly implies the continuity of $K_{2,p}$.
We next turn to the proof of the continuity of $K_{m,p}$ for odd $m$. 
Consider a polynomial $\pi\in \H_m$. If $K_{m,p}(\pi)=0$ then the upper semi-continuity of $K_{m,p}$, combined with its non-negativity, implies that it is continuous at $\pi$. 
Otherwise, assume that $K_{m,p}(\pi)>0$. Consider a sequence $\pi_n\in \H_m$ converging to $\pi$, and a sequence $\phi_n$ of linear transformations satisfying $\det \phi_n = 1$, and such that
$$
\lim_{n\to +\infty} e_{\phi_n(\Teq)}(\pi_n) =  \liminf_{\pi^*\to \pi} K_{m,p}(\pi^*)
:= \lim_{r\to 0} \inf_{\|\pi- \pi\|\leq r} K_{m,p}(\pi^*).
$$
If the sequence $\phi_n$ admits a converging subsequence $\phi_{n_k}\to \phi$, it follows that 
$$
K_{m,p}(\pi) \leq e_{\phi(\Teq)}(\pi) = \lim_{k\to+\infty} e_{\phi_{n_k}(\Teq)}(\pi_{n_k}) = \liminf_{\pi^*\to \pi} K_{m,p}(\pi^*).
$$
This asserts that $K_{m,p}$ is lower semi continuous at $\pi$, and therefore continuous at $\pi$ since we already know that $K_{m,p}$ is upper semi-continuous.

If $\phi_n$ does not admit any converging subsequence, then we invoke the SVD decomposition 
$\phi_n = U_n\circ D_n \circ V_n$, where $U_n,V_n\in \cO_2$ and $D_n = \diag(\ve_n,\frac 1 {\ve_n})$, where $0<\e_n\leq 1$. (Here and below, we use the shorthand $\diag(a,b)$ to denote the diagonal matrix with entries $a$ and $b$) 
The compactness of $\cO_2$ implies that $U_n$ admits a converging subsequence $U_{n_k}\to U$. 
In particular $\pi_{n_k}\circ U_{n_k}$ converges to $\pi\circ U$. Therefore,
denoting by $a_{i,n}$ the coefficient of $x^i y^{m-i}$ in $\pi_n \circ U_n$, the subsequence $a_{i,n_k}$ converges to the coefficient $a_i$ of $x^i y^{m-i}$ in $\pi\circ U$.
Observe also that $\e_n\to 0$, otherwise some converging 
subsequence could be extracted from $\phi_n$.
Since $e_{\phi_n(\Teq)}(\pi_n) = e_{\Teq}(\pi_n\circ\phi_n)$, the sequence of polynomials $\pi_n\circ\phi_n$ is uniformly bounded, and so is the sequence $\pi_n\circ U_n\circ D_n$. Therefore
the sequences $(a_{i,n} \ve_n^{2i-m})_{n\geq 0}$ are uniformly bounded.
It follows that $a_i = 0$ when $i< \frac m 2$. Since $m$ is odd, this implies that $\pi\circ U(x,y) = x^{s_m} \tilde\pi(x,y)$ and Proposition \ref{vanishprop} implies that $K_{m,p}(\pi) = 0$ which contradicts the hypothesis $K_{m,p}(\pi)>0$.

Last, we prove property \iref{contzero}.
The assumption $K_{m,p}(\pi_n)\to 0$ is equivalent to the existence of 
a sequence $T_n=\phi_n(\Teq)$ with $\det \phi_n = 1$ such that
$e_{m,T_n}(\pi_n)_p \to 0$. Reasoning in a similar way as in
the proof of Proposition \ref{vanishprop}, we first obtain that $\pi_n\circ \phi_n\to 0$,
and we then invoke the SVD decomposition of $\phi_n$
to build a converging sequence of orthogonal matrices 
$U_n\ra U$ and a sequence $0<\ve_n \leq 1$ such that if $a_{i,n}$ is the coefficient of 
$x^i y^{m-i}$ in $\pi_n\circ U_n$, we have $a_{i,n} \ve_n^{2i-m} \to 0$. 
When $i<s_m$, it follows that $a_{i,n}\to 0$ and therefore $\pi\circ U (x,y)= x^{s_m} \hat \pi(x,y)$. The result follows from Proposition \ref{vanishprop}.
\sq

We finally make a connection between the shape function
and the approach developed in \cite{C3}.
For all $\pi\in \H_m$, we denote by $\Lambda_\pi$ as the level set of $|\pi|$ for the value $1$,
\be
\Lambda_\pi = \{(x,y)\in \R^2,\ |\pi (x,y)|\leq 1\}.
\label{lambdapi}
\ee
We now define 
\be
\label{defKE}
K^\cE_m(\pi) = \left(\sup_{E\in \cE,\ E\subset \Lambda_\pi} |E|\right)^{-m/2}.
\ee
where the supremum is taken over the set $\cE$ of all ellipses centered at $0$.
The optimization problem defining $K^\cE_m$ is equivalent to 
\be
\label{OptimEll}
 \inf \{\det H\sep  H\in S_2^+ \text{ and } \forall z \in \R^2, \<Hz,z\> \geq |\pi(z)|^{2/m}\},
\ee
where $S_2^+$ is the cone of $2\times 2$ symmetric definite positive matrices.
The minimizing ellipse $E^*$ is then given by $\{ \<Hz,z\>\leq 1\}$.
The optimization problem described in \iref{OptimEll} 
is quadratic in dimension $2$, and subject to (infinitely many) linear constraints.
This apparent simplicity is counterbalanced by the fact that it is non convex. In
particular, it does not have unique solutions and may also have no solution.

\begin{prop}
\label{propequivEllTri}
On $\H_m$, one has the equivalence

$$
cK^\cE_m \leq K_{m,p} \leq C K^\cE_m
$$
with constant $0<c\leq C$ independent of $p$.
\end{prop}

\proof Let $\Teq$ denote an equilateral triangle of unit area, 
and $B$ its circumscribed disk. It is easy to see that the inscribed disc is $B/2$.

We first show that $K_{m,p} \leq CK_{m}^{\cE}$.
Let $\pi\in \H_m$, and let $E_n$ be a sequence of ellipsoids inscribed in $\Lambda_\pi$
and such that $|E_n|$ tends to $\sup \{|E|\sep E\in \cE,\ E\subset \Lambda_\pi\}$  as $n\to +\infty$. 
We write $E_n = \lambda_n \phi_n(B)$, where $\phi_n$ is a linear transform such that
$\det \phi_n=1$  and $\lambda_n>0$. We define the triangle $T_n = \phi_n(\Teq)$ which 
satisfies $|T_n|=1$. We then have
\begin{eqnarray*}
K_{m,p}(\pi)  &\leq &
\|\pi - I_{m,T_n} \pi\|_{L^p(T_n)} \\
&= & \|\pi\circ\phi_n - (I_{m,T_n} \pi)\circ\phi_n\|_{L^p(\Teq)}\\
&= & \|\pi\circ\phi_n - I_{m,\Teq} (\pi\circ\phi_n)\|_{L^p(\Teq)}\\
&\leq & \|\pi\circ\phi_n - I_{m,\Teq} (\pi\circ\phi_n)\|_{L^\infty(\Teq)},
\end{eqnarray*}
where we have used the commutation formula \iref{intercommut}.

Remarking that $I_{m,\Teq}$ is a continuous
operator from $\H_m$ to $\P_{m-1}$ in the sense of any norm since these
spaces are finite dimensional, we thus obtain
\begin{eqnarray*}
K_{m,p}(\pi)  
& \leq & C_1 \|\pi\circ\phi_n\|_{L^\infty(\Teq)}\\
& \leq & C_1 \|\pi\circ\phi_n\|_{L^\infty(B)}\\
&=& C_1 \|\pi\|_{L^\infty(\phi_n(B))}\\
&=& C_1 \lambda_n^{-m} \|\pi\|_{L^\infty(E_n)}\\
&\leq & C_1 \left(\frac{|E_n|}{|B|}\right)^{-m/2},
\end{eqnarray*}
where we have used the fact that $|\pi|\leq 1$ in $E_n\subset \Lambda_\pi$.
Letting $n\to +\infty$, we obtain that $K_{m,p}(\pi)\leq C K_m^\cE(\pi)$ with
$C= C_1 |B|^{m/2}$.

We next prove that $cK_m^\cE \leq K_{m,p}$.
Let $T_n$ be a sequence of triangles of unit area such that
$e_{m,T_n}(\pi)_p$ tends to $K_{m,p}(\pi)$ as $n\to +\infty$.
As already remarked in \iref{transinv}
the interpolation error is invariant by translation.
We may therefore assume that the triangles $T_n$ have their barycenter at the origin. 
Then there exists linear transforms $\phi_n$ with $\det \phi_n=1$, such that $T_n = \phi_n(\Teq)$. 
We now write
\begin{eqnarray*}
\|\pi\|_{L^\infty(\phi_n(B/2))} &\leq & \|\pi\|_{L^\infty(T_n)}\\
&=& \|\pi\circ\phi_n\|_{L^\infty(\Teq)}\\
&\leq & C_2 e_{m,\Teq}(\pi\circ\phi)_1,\\
&\leq & C_2 e_{m,\Teq}(\pi\circ\phi)_p,\\ 
& = & C_2  e_{m,T_n}(\pi)_p,
\end{eqnarray*} 
where we have used the fact that $\|\cdot\|_{L^\infty(\Teq)}$
and $e_{\Teq}(\cdot)_1$ are equivalent norms on $\H_m$, and that  $e_{m,\Teq}(\cdot)_p$ is an increasing function of $p$ since $|\Teq| = 1$. By homogeneity, it 
follows that if  $\lambda_n:= (C_2 e_{m,T_n}(\pi)_p)^{-1/m}$, we have
$$
\|\pi\|_{L^\infty(\lambda_n \phi_n(B/2))} \leq  \lambda_n^m C_2  e_{m,T_n}(\pi)_p=1.
$$
Therefore the ellipse $E_n:=\lambda_n \phi_n(B/2)$
is contained in $\Lambda_\pi$, so that
$$
K_m^\cE \leq |E_n|^{-m/2} =\left(\lambda_n^2 \frac {|B|} 4\right)^{-m/2}=C_2\left(\frac4 {|B|} \right)^{m/2} e_{m,T_n}(\pi)_p.
$$
Letting $n\to +\infty$, we obtain that $cK_m^\cE \leq K_{m,p}$
with $c=C_2^{-1}\left(\frac {|B|} 4\right)^{m/2}$.
\sq

\begin{remark}
Since $K_{m,p}$ and $K_m^\cE$ are equivalent, they must vanish on the
same set, and therefore Proposition \ref{vanishprop} is also valid for $K_m^\cE$.
It also easy to see that $K_m^\cE$ satisfies the homogeneity and invariance
properties stated for $K_{m,p}$ in \iref{homogK} and \iref{invK}, as well as
the continuity properties stated in Proposition \ref{propsemicont}.
\end{remark}

\begin{remark}
The continuity of the functions $K_{m,p}$ and $K_m^\cE$ can be established
when $m$ is odd or equal to $2$, as shown by Proposition \ref{propsemicont}, but seems to fail otherwise. In particular, direct computation shows that $K_4^\cE(x^2y^2-\ve y^4)$ is independent of 
$\ve>0$ and strictly smaller than $K_4^\cE(x^2 y^2)$. 
Therefore $K_4^\cE$ is upper semi-continuous but discontinuous at the point $x^2 y^2\in \H_4$.
\end{remark}

\section{Optimal estimates}

This section is devoted to the proofs of our main theorems, starting with the lower estimate of Theorem \ref{optitheorem}, and continuing with the upper estimates involved in both Theorem \ref{maintheorem} and \ref{optitheorem}.

Throughout this section, for the sake of notational simplicity, 
we fix the parameters $m$ and $p$ and use the shorthand
$$
K=K_{m,p}\;\; {\rm and}\;\; e_T(\pi) = e_{m,T}(\pi)_p.
$$
For each point $z\in \Omega$ we define
$$
\pi_z := \frac{d^mf_z}{m!} \in \H_m,
$$
where $f\in C^m(\Omega)$ is the function in the statement of the theorems. We 
denote by
$$
\omega(r):=\sup_{\|z-z'\|\leq r} \|\pi_z-\pi_{z'}\|,
$$
the modulus of continuity of $z\mapsto \pi_z$
with the norm $\|\cdot\|$ defined by \iref{normdef}. 
Note that $\omega(r)\to 0$ as $r\to 0$.

\subsection{Lower estimate}

In this proof we will use an estimate by below of the local interpolation error.

\begin{prop}
\label{errorfPz}
Assume that $1\leq p<\infty$. There exists a constant $C>0$, depending on $f$ and $\Omega$,
such that for all triangle $T\subset \Omega$ and $z\in T$,
\be
e_T(f)^p \geq K^p(\pi_z) |T|^{\frac{mp} 2+1} - C(\diam T)^{mp} |T| \omega(\diam T).
\label{localbelow}
\ee
\end{prop}

\proof
Denoting by $\mu_z\in \P_m$ the Taylor development of $f$ at the point $z$ up to degree $m$, we obtain
$$
f(z+u) - \mu_z(z+u) = m \int_{t=0}^1 (\pi_{z+tu}(u)-\pi_z(u))(1-t)^{m-1} dt.
$$
and therefore
$$
\|f-\mu_z\|_{L^\infty(T)} \leq C_0\diam(T)^m \omega(\diam(T)),
$$
where $C_0$ is a fixed constant.
By construction $\pi_z$ is the homogenous part of $\mu_z$ of degree $m$, and therefore $\mu_z-\pi_z\in \P_{m-1}$. It follows that for any triangle $T$, we have 
\be
\label{eqMuPi}
\mu_z - I_{m,T}\mu_z = \pi_z -I_{m,T} \pi_z.
\ee
We therefore obtain 
\begin{eqnarray*}
|e_T(f)- e_T(\pi_z)| &\leq & \|(f-I_{m,T}f) - (\pi_z-I_{m,T}\pi_z)\|_{L^p(T)} \\
&\leq & |T|^{1/p} \|(f-I_{m,T}f) - (\mu_z-I_{m,T}\mu_z)\|_{L^\infty(T)} \\
&= & |T|^{1/p} \|(I-I_{m,T})(f-\mu_z)\|_{L^\infty(T)} \\
&\leq & C_1|T|^{1/p} \|f-\mu_z\|_{L^\infty(T)} \\
&\leq & C_0C_1 |T|^{1/p} \diam(T)^m \omega(\diam(T))
\end{eqnarray*}
where $C_1$ is the norm of the operator $I-I_{m,T}$ in $L^\infty(T)$
which is independent of $T$. 

From \iref{shapeA} we know that $e_T(\pi_z) \geq |T|^{\frac m 2+\frac 1 p} K(\pi_z)$, and therefore
$$
e_T(f) \geq K(\pi_z) |T|^{\frac m 2+\frac 1 p} - C_0C_1 |T|^{1/p} \diam(T)^m \omega(\diam(T)).
$$
We now remark that for all $p\in[1,\infty)$ the function $r\mapsto r^p$ is convex, and therefore if $a,b,c$ are positive numbers, and $a\geq b-c$ then $a^p \geq \max\{0,b-c\}^p \geq b^p -p c b^{p-1}$. Applying this to our last inequality we obtain
$$
e_T(f)^p \geq K^p(\pi_z) |T|^{\frac {mp} 2+1} -
pC_0C_1 (K(\pi_z))^{p-1} |T|^{(p-1)(\frac m 2+\frac 1 p)+\frac 1 p} \diam(T)^m \omega(\diam T).
$$
Since $|T|^{(p-1)(\frac m 2+\frac 1 p)+\frac 1 p}=|T|^{(p-1)\frac m 2} |T| \leq (\diam T)^{m(p-1)}|T|$, this leads to
$$
e_T(f)^p \geq K^p(\pi_z) |T|^{\frac{mp} 2+1} - C(\diam T)^{mp} |T| \omega(\diam T),
$$
where $C :=pC_0C_1 (\sup_{z\in\Omega}K(\pi_z))^{p-1}$. 
\sq

We now turn to the proof of the lower estimate in Theorem \ref{optitheorem} in the case
where $p<\infty$.
Consider a sequence $\seqT$ of triangulations which is admissible in the sense of equation 
\iref{admissibilitycond}. Therefore, there exists a constant $C_A$ such that 
$$
\diam T\leq C_A N^{-1/2},\; N\geq N_0,\; T\in \cT_N
$$
For $T\in \cT_N$, we combine this estimate with \iref{localbelow}, which gives
$$
e_T(f)^p \geq K^p(\pi_z) |T|^{\frac{mp} 2+1} - (C_AN^{-1/2})^{mp} |T| C\omega(C_AN^{-1/2}).
$$
Averaging over $T$, we obtain
$$
e_T(f)^p \geq \int_T K^p(\pi_z) |T|^{\frac{mp} 2} dz - |T|  N^{-\frac{mp} 2}C_A^{mp}
C\omega(C_A N^{-1/2}).
$$
Summing on all $T\in \cT_N$, and denoting by $T_z^N$ the triangle in $\cT_N$ containing the point $z\in\Omega$, we obtain the estimate
\be
\label{lowerint}
e_{\cT_N}(f)^p \geq \int_\Omega K(\pi_z) |T_z^N|^{\frac{mp} 2} dz - N^{-\frac{mp} 2}\e(N),
\ee
where $\e(N):=|\Omega|C_A^{mp}C\omega(C_AN^{-1/2})\to 0$ as $N\to +\infty$.
The function $z\mapsto |T_z^N|$ is linked with the number of triangles in the following way:
$$
\int_\Omega \frac {dz} {|T_z^N|} = \sum_{T\in \cT_N} \int_T \frac 1{|T|} = N.
$$
On the other hand, with $\frac 1 q = \frac m 2 +\frac 1 p$, we have by H\"older's inequality,
\be
\int_\Omega K^q(\pi_z) dz\leq \left(\int_\Omega K^p(\pi_z) |T_z^N|^{\frac{mp} 2} dz\right)^{q/p}
\left(\int_\Omega \frac 1 {|T_z^N|}dz\right)^{1-q/p}.
\label{holder}
\ee
Combining the above, we obtain a lower bound for the integral term in \iref{lowerint}
which is independent of $\cT_N$:
$$
\int_\Omega K^p(\pi_z) |T_z^N|^{\frac{mp} 2} dz\geq \left(\int_\Omega K^q(\pi_z)dz\right)^{p/q} N^{-m p/2}.
$$
Injecting this lower bound in \iref{lowerint} we obtain
$e_{\cT_N}(f)^p \geq \left[\left(\int_\Omega K^q(\pi_z)dz\right)^{p/q} -\e(N) \right] N^{-m p/2}.$
This allows us to conclude
\be
\liminf_{N\to +\infty} N^{\frac m 2} e_{\cT_N}(f) \geq \left(\int_\Omega K^q(\pi_z)dz\right)^{\frac 1 q},
\label{lowerest}
\ee
which is the desired estimate.
\nl
\nl
The case $p=\infty$ follows the same ideas. 
Adapting Proposition \ref{errorfPz}, one proves that
$$
e_T(f) \geq K(\pi_z) |T|^{\frac m 2} - C(\diam T)^{m} \omega(\diam T).
$$
and therefore
\be
e_{\cT_N}(f) \geq \left\|K(\pi_z) |T_z^N|^{\frac m 2} \right\|_{L^\infty(\Omega)} - N^{-\frac m 2} \e(N), 
\label{lowerint2}
\ee
where $\e(N):= C_A^m C\omega(C_A N^{-\frac 1 2})\to 0$ as $N\to +\infty$.
The Holder inequality now reads
$$
\int_\Omega K(\pi_z)^{\frac 2 m} dz \leq \left\|K(\pi_z)^{\frac 2 m} |T_z^N|\right\|_{L^\infty(\Omega)} \left\|\frac 1 {|T_z^N|}\right\|_{L^1(\Omega)} 
$$
equivalently
$$
\left\|K(\pi_z) |T_z^N|^{\frac m 2} \right\|_{L^\infty(\Omega)} \geq \left(\int_\Omega K(\pi_z)^{\frac 2 m} dz\right)^{\frac m 2} N^{- \frac m 2}.
$$
Combining this with \iref{lowerint2}, this leads to the desired estimate \iref{lowerest}
with $p=\infty$ and $q=\frac 2 m$.

\begin{remark}
This proof reveals the two principles which characterize the optimal triangulations.
Indeed, the lower estimate \iref{lowerest}
becomes an equality only when both inequalities in \iref{localbelow} and \iref{holder}
are equality. The first condition - equality in \iref{localbelow} - is met when each 
triangle $T$ has an optimal shape, in the sense that $e_T(\pi_z)=K(\pi_z)|T|^{\frac m 2+\frac 1 p}$
for some $z\in T$. The second condition - equality in \iref{holder} - is met when the ratio
between 
$K^p(\pi_z)|T_z^N|^{\frac{mp} 2}$ and $|T_z^N|^{-1}$ is constant, or equivalently
$K(\pi_z)|T|^{\frac m 2+\frac 1 p}$ is independent of the triangle $T$. 
Combined with the first condition, this means that the error $e_T(f)^p$ is equidistributed
over the triangles, up to the perturbation by  $(\diam T)^{mp} |T| \omega(\diam T)$ which
becomes neglectible as $N$ grows.
\end{remark}

\subsection{Upper estimate}

We first remark that the upper estimate in Theorem \ref{optitheorem}. 
implies the upper estimate in Theorem \ref{maintheorem}
by a sub-sequence extraction argument: if the upper estimate in Theorem \ref{optitheorem}
holds, then for all $n>0$ there exists a sequence
$(\cT_N^n)_{N>N_0}$ such that 
$$
\limsup_{N\to +\infty}\( N^{\frac m 2} e_{\cT_N^n}(f)\) \leq \left\| K\left(\frac {d^mf}{m!}\right)\right\|_{L^q}+ \frac 1 n,
$$
with $\frac 1 q=\frac 1 p+\frac m 2$. We then take $\cT_N=\cT_N^{n(N)}$, where 
$$
n(N)=\max\left\{n\leq N\; ;\;  N^{\frac m 2} e_{\cT_N^n}(f)\leq \left\| K\left(\frac {d^mf}{m!}\right)\right\|_{L^q}+ \frac 2 n\right\}.
$$
For $N$ large enough this set is finite and non empty, and therefore $n(N)$ is well defined. Furthermore $n(N)\to +\infty$ as $N\to +\infty$ and therefore 
$$
\limsup_{N\to +\infty}\( N^{\frac m 2} e_{\cT_N}(f)\) \leq \left\| K\left(\frac {d^mf}{m!}\right)\right\|_{L^q}.
$$
We are thus left with proving the upper estimate in Theorem \ref{optitheorem}. 
We begin by fixing a (large) number $M>0$. We shall take the limit $M\to \infty$ in the very last step of our proof. We define 
$$
\mT_M = \{T \text{ triangle, } |T|=1,\ \bary(T)=0 \text{ and } \diam(T)\leq M\},
$$
the set of triangles centered at the origin, of unit area and diameter smaller than $M$.
This set is compact with respect to the Hausdorff distance.
This allows us to define a ``tempered'' version of $K = K_{m,p}$ that we denote by $K_M$:
$$
K_M(\pi) = \inf_{T\in \mT_M} e_T(\pi).
$$
Since $\mT_M$ is compact, the above infimum
is attained on a triangle that we denote by $T_M(\pi)$.
Note that the map $\pi\mapsto T_M(\pi)$ need not be continuous.
It is clear that $K_M(\pi)$ decreases as $M$ grows.
Note also that the restriction to triangles $T$ centered at $0$ is artificial, since
the error is invariant by translation as noticed in 
\iref{transinv}. Therefore $K_M(\pi)$ converges to $K(\pi)$ as $M\to+\infty$.
Since $\mT_M$ is compact, the map $\pi \mapsto \max_{T\in\mT_M}e_T(\pi)$ 
defines a norm on $\H_m$, and is therefore bounded by $C_M\|\pi\|$ for some $C_M>0$.
One easily sees that the functions $\pi \mapsto e_T(\pi)$ are uniformly $C_M$-Lipschitz
for all $T\in \mT_M$, and so is $K_M$.

We now use this new function $K_M$ to obtain a local upper error estimate that is closely related to the local lower estimate in Proposition \ref{errorfPz}
\begin{prop}
For $z_1\in\Omega$, let $T$ be a triangle which is obtained from
$T_M(\pi_{z_1})$ by rescaling and translation ( $T=tT_M(\pi_{z_1})+z_0$).
Then for any $z_2\in T$, 
\be
\label{localabove}
e_T(f) \leq \(K_M(\pi_{z_2}) + B_M\omega(\max\{|z_1-z_2|,\diam(T)\})\) |T|^{\frac m 2+\frac 1 p},
\ee
where $B_M>0$ is a constant which depends on $M$.
\end{prop}

\proof

For all $z_1,z_2\in \Omega$, we have
\begin{eqnarray*}
\label{eTMPxPy}
e_{T_M(\pi_{z_1})}(\pi_{z_2}) &\leq & e_{T_M(\pi_{z_1})} (\pi_{z_1}) + C_M \|\pi_{z_1}-\pi_{z_2}\| \\
& = & K_M(\pi_{z_1}) + C_M \|\pi_{z_1}-\pi_{z_2}\|,\\
&\leq & K_M(\pi_{z_2}) + 2 C_M\|\pi_{z_1}-\pi_{z_2}\|,\\
&\leq & K_M(\pi_{z_2}) + 2C_M\omega (|z_1-z_2|).
\end{eqnarray*}
Therefore, if $T$ is of the form $T=tT_M(\pi_{z_1})+z_0$, we obtain by a change of variable that
$$
e_T(\pi_{z_2})\leq \(K_M(\pi_{z_2}) + 2C_M\omega(|z_1-z_2|)\) |T|^{\frac m 2+\frac 1 p} 
$$
Let $\mu_z\in \P_m$ be the Taylor polynomial of $f$ at the point $z$ up to degree $m$.
Using Equation \iref{eqMuPi} we obtain 
\begin{eqnarray*}
e_T(f) &\leq &  e_T(\mu_{z_2})+ e_T(f-\mu_{z_2})\\
&=& e_T(\pi_{z_2})+ e_T(f-\mu_{z_2})\\
&\leq &  \(K_M(\pi_{z_2}) + 2C_M\omega(|z_1-z_2|)\) |T|^{\frac m 2+\frac 1 p} + e_T(f-\mu_{z_2})\\
\end{eqnarray*}
By the same argument as in the proof of Proposition \ref{errorfPz}, we derive that
$$
e_T(f-\mu_{z_2})\leq C|T|^{\frac 1 p} \diam(T)^m \omega(\diam T),
$$
and thus 
$$
e_T(f) \leq \(K_M(\pi_{z_2}) + 2C_M\omega(|z_1-z_2|)\) |T|^{\frac m 2+\frac 1 p} +
C|T|^{\frac 1 p} \diam(T)^m \omega(\diam T).
$$
Since $T$ is the scaled version of a triangle in $\mT_M$, it obeys $\diam(T)^2 \leq M^2 |T|$. Therefore 
$$
e_T(f) \leq \left(K_M(\pi_{z_2}) + (2C_M+CM^m)\omega(\max\{|z_1-z_2|,\diam (T)\})\right) |T|^{\frac m 2 +\frac 1 p},
$$
which is the desired inequality with $B_M:=2C_M+CM^m$.
\sq

For some $r>0$ to be specified later, we now choose an arbitrary triangular mesh $\cR$ of $\Omega$ satisfying 
$$
r \geq \sup_{R\in \cR} \diam(R).
$$
Our strategy to build a triangulation that satisfies the optimal upper estimate
is to use the triangles $R$ as {\it macro-elements} in the sense that each of them
will be tiled by a locally optimal uniform triangulation.  This strategy was already used in \cite{B}.

For all $R\in \cR$ we consider the triangle
$$
T_R:=(K_M(\pi_{b_R})+2B_M\omega(r))^{-\frac q 2} T_M(\pi_{b_R}),
$$
which is a scaled version of $T_M(\pi_{b_R})$
where $b_R$ is the barycenter of $R$. We use this triangle
to build a periodic tiling $\cP_R$
of the plane: there exists a vector $c$ such that 
$T_R\cup T'_R$ forms a parallelogram
of side vectors $a$ and $b$, with $T'_R=c-T_R$. We
then define 
\be
\cP_R:=\{T_R+ma+nb\sep m,n\in\Z^2\} \cup  \{T'_R +ma+nb\sep m,n\in\Z^2\}.
\label{defTiling}
\ee
Observe that for all $\pi\in \H_m$, and all triangles $T,T'$ such that $T'=-T$ one has $e_T(\pi) = e_{T'}(\pi)$ since $\pi$ is either an even polynomial when $m$ is an even integer, or an odd polynomial when $m$ is odd. Since we already know that $e_T(\pi)$ is invariant by translation of $T$,
we find that the local error $e_T(\pi)$ is constant on all $T\in \cP_R$.

We now define as follows a family of triangulations $\cT_s$ of the domain 
$\Omega$, for $s>0$. For every $R\in \cR$, we consider the elements $T\cap R$
for $T\in s\cP_R$, where $s\cP_R$ denotes the triangulation $\cP_R$
scaled by the factor $s$. Clearly $\{T\cap R,\; T\in s\cP_R,\; R\in\cR\}$ constitute
a partition of $\Omega$. In this partition, we distinguish the interior 
elements 
$$
\cTr:=\{T\in s\cP_R \sep T\in {\rm int}(R) \; , R\in \cR\},
$$
which define pieces of a conforming triangulation, and the 
boundary elements $T\cap R$ for $T\in s\cP_R$ such that $T\cap \partial R\neq \emptyset$.
These last elements might not be triangular, nor conformal with the 
elements on the other side. Note that for $s>0$ small enough, each $R\in\cR$
contains at least one triangle in $\cTr$, and therefore the boundary elements
constitute a layer around the edges of $\cR$. In order to obtain a conforming
triangulation, we proceed as follow: for each boundary element $T\cap R$,
we consider the points on its boundary which are either its vertices or
those of a neighboring element. We then build the Delaunay triangulation
of these points, which is a triangulation of $T\cap R$ since it is a convex set.
We denote by $\cTb$ the set of all triangles obtained by this procedure,
which is illustrated on Figure 1.

\begin{figure}
	\centering
		\includegraphics[width=4cm,height=4cm]{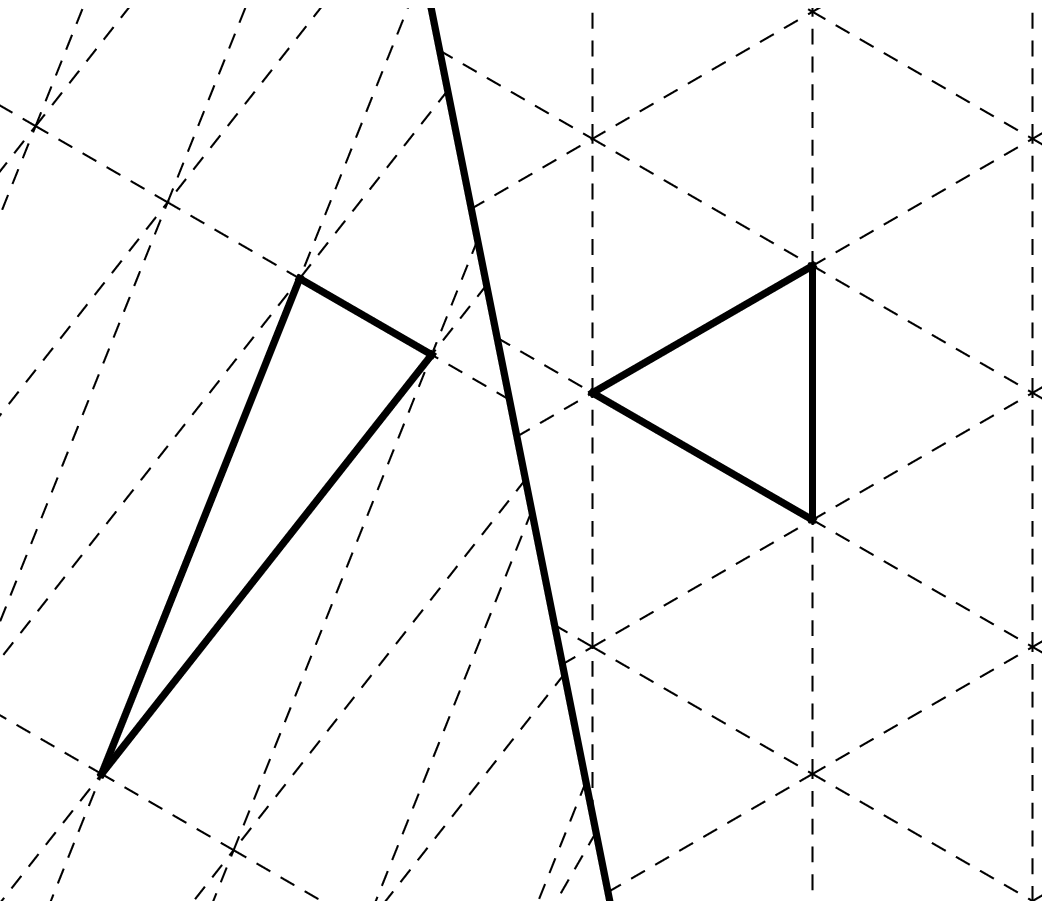}
		\hspace{1cm}
		\includegraphics[width=4cm,height=4cm]{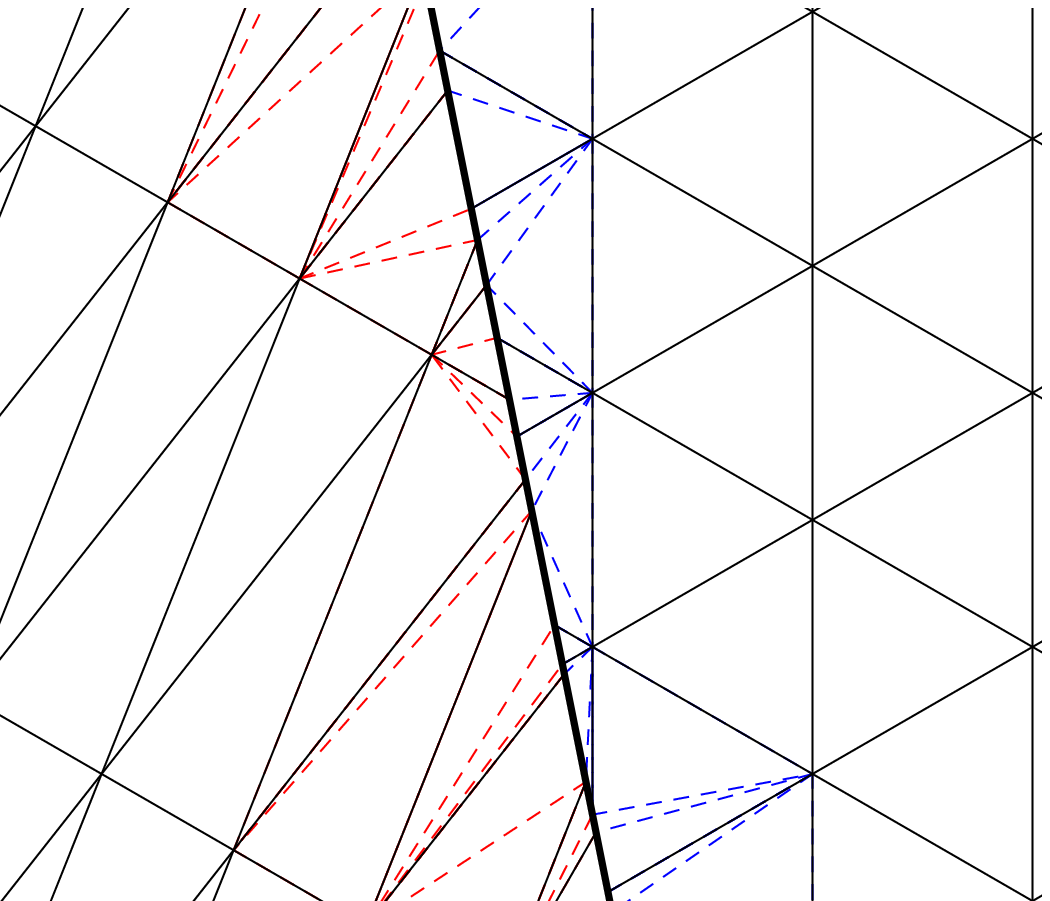}
		\hspace{1cm}
		\includegraphics[width=4cm,height=4cm]{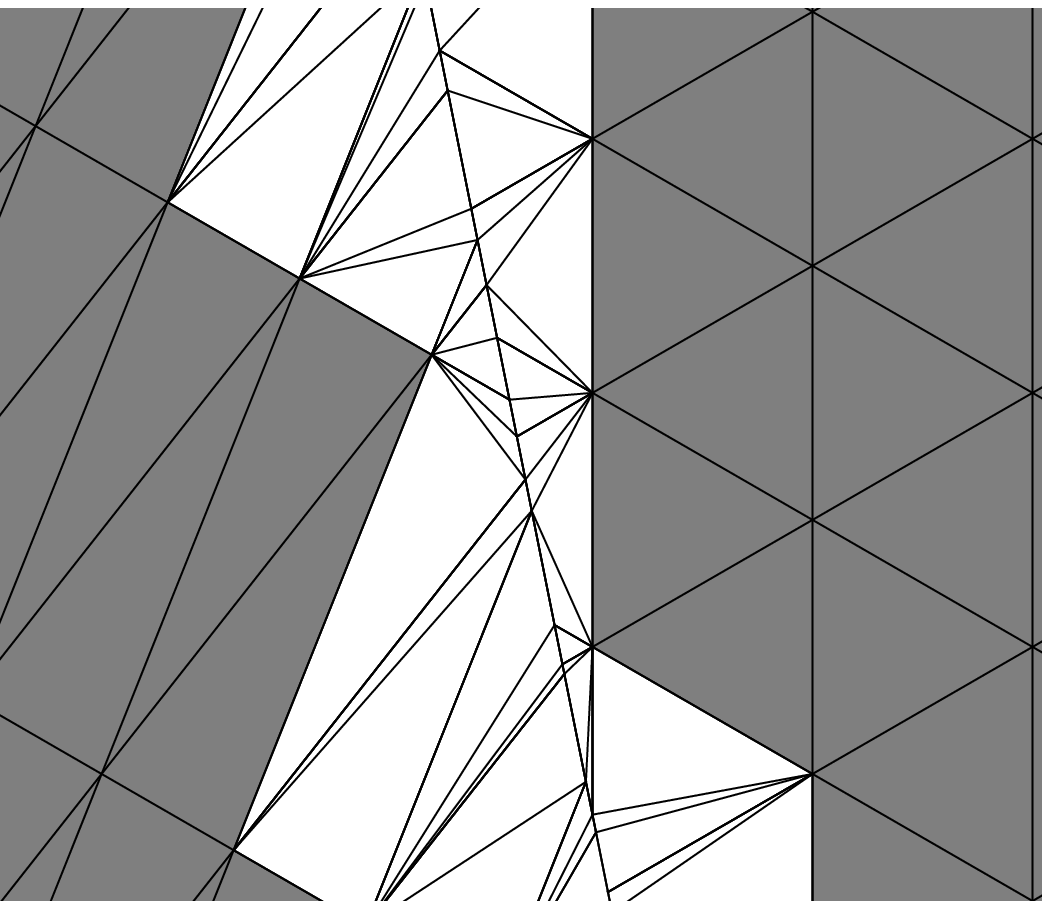}
	\caption{a. An edge (Thick) of the macro-triangulation $\cR$ separating to uniformly paved regions ($T_R$ is thick, $\cP_R$ is dashed). b. Additional edges (dashed) are added near the interface in order to preserve conformity. c. The sets of triangles $\cTr$ (gray) and $\cTb$ (white)}
\end{figure}

Our conforming triangulation is given by
$$
\cT_s=\cTr \cup \cTb.
$$
As $s\ra 0$, clearly 
$$
\#(\cTb)\leq  C_{\rm bd}s^{-1} \text{ and } \sum_{T\in \cTb} |T| \leq C_{\rm bd} s,
$$
for some constant $C_{\rm bd}$ which depends on the macro-triangulation $\cR$. 
We do not need to estimate $C_{\rm bd}$ since $\cR$ is fixed and the contribution due to $C_{\rm bd}$ in the following estimates is neglectible as $s\to 0$.
We therefore obtain that the number of triangles
in $\cTb$ is dominated by the number of triangles in $\cTr$. More precisely, we
have the equivalence
\be
\label{cardcTs}
\#( \cT_s )\sim \#(\cTr) \sim  \sum_{R\in \cR} \frac{ |R|}{s^2|T_R|}
= s^{-2} \sum_{R\in \cR} |R| (K_M(\pi_{b_R})+2B_M \omega(r))^q,
\ee
in the sense that the ratio between the above quantities tends to $1$ as $s\to 0$.
The right hand side in \iref{cardcTs} can be estimated through an integral:
\begin{eqnarray*}
s^2 \#(\cTr) &\leq & \sum_{R\in \cR} |R| (K_M(\pi_{b_R})+2B_M \omega(r))^q\\
& = & \sum_{R\in \cR} \int_R (K_M(\pi_{b_R})+2B_M \omega(r))^q dz\\
& \leq & \sum_{R\in \cR} \int_R (K_M(\pi_z)+C_M\|\pi_z-\pi_{b_R}\|+2B_M \omega(r))^q dz\\
&\leq & \int_\Omega (K_M(\pi_z)+(2B_M+C_M) \omega(r))^q dz
\end{eqnarray*}
Therefore, since $C_M\leq B_M$,
\be
\label{cardcTsUpper}
\#(\cT_s)\leq s^{-2}\left(\int_\Omega (K_M(\pi_z)+3B_M \omega(r))^q dz +C_{\rm bd} s \right).
\ee
Observe that the construction of $\cT_s$ gives a bound on the diameter of its elements
$$
\sup_{T\in \cT_s}\diam(T)\leq s C_a, \;\; C_a:=\max_{R\in\cR}{\rm diam}(T_R)
$$
Combining this with \iref{cardcTs}, we obtain that   
$$
\sup_{T\in \cT_s}\diam(T) \leq C_A (\#\cT_s)^{-1/2} \text{ for all } s>0
$$
which is analogous to the admissibility condition 
\iref{admissibilitycond}. 
\nl
\nl
We now estimate the global interpolation error 
$\|f-I_{m,\cT_s}f\|_{L^p}:=(\sum_{T\in\cT_s}e_T(f)^p)^{\frac 1 p}$,
assuming first that $1\leq p<\infty$. We first estimate
the contribution of $\cTb$, which will eventually be neglectible. 
Denoting $\nu_z\in \P_{m-1}$ the Taylor polynomial of $f$ up to degree $m-1$ at $z$ we remark that
$$
\|f-I_{m,T} f\|_{L^\infty(T)} = \|(I-I_{m,T} )(f-\nu_{b_T})\|_{L^\infty(T)}
 \leq C_1\|f-\nu_{b_T}\|_{L^\infty(T)} \leq C_0C_1 \diam(T)^m.
$$
where $C_1$ is the norm of $I-I_{m,T}$
in $L^\infty(T)$ which is independent of $T$
and $C_0$ only depends on the $L^\infty$ norm
of $d^mf$. Remarking that $e_T(f) = \|f-I_{m,T} f\|_{L^p(T)} \leq |T|^{\frac 1 p} \|f-I_{m,T} f\|_{L^\infty(T)}$, we obtain an upper bound for the contribution of $\cTb$ to the error:
\begin{eqnarray*}
\sum_{T\in \cTb}  e_{T}(f)^p &\leq& C_0^pC_1^p\sum_{T\in \cTb} |T|  \diam(T)^{mp}\\
&\leq & C_0^pC_1^p \left(\sum_{T\in \cTb} |T|\right) \sup_{T\in \cTb} \diam(T)^{mp} \\
&\leq & C_0^pC_1^p C_{\rm bd} s \sup_{T\in \cTb} \diam(T)^{mp}\\
&\leq & C^*_{\rm bd}s^{mp+1},
\end{eqnarray*}
with $C^*_{\rm bd}=C_0^pC_1^p C_a^{mp}C_{\rm bd}$.
We next turn to the the contribution of $\cTr$ to the error.
If $T\in \cTr$, $T\subset R\in \cR$, we consider any point $z_1 = z\in T$ and define $z_2 = b_R$
the barycenter of $R$. With such choices, the estimate \iref{localabove} reads 
$$
e_{T}(f) \leq \left(K_M(\pi_z) + B_M \omega(\max\{r,C_A s\})\right) |T|^{\frac m 2+\frac 1 p}.
$$
We now assume that $s$ is chosen small enough such that 
$C_A s\leq r$. Geometrically, this condition ensures that the ``micro-triangles'' constituting $\cT_s$ actually have a smaller diameter than the ``macro-triangles'' constituting $\cR$. This implies
\be
\label{estimabove}
e_{T}(f)^p \leq \left(K_M(\pi_z) + B_M \omega(r)\right)^p |T|^{\frac {mp} 2+1}
\ee
Given a triangle $T\in \cTr$, $T\subset R \in \cR$, and a point $z\in T$, one has 
\begin{eqnarray*}
|T|  &= &s^2 \left(K_M(\pi_{b_R})+2B_M \omega(r)\right)^{-q}\\
& \leq & s^2 \left(K_M(\pi_{z})-C_M\|\pi_z-\pi_{b_R}\|+2B_M \omega(r)\right)^{-q}\\
&\leq& s^2\left(K_M(\pi_z)+(2B_M-C_M) \omega(r)\right)^{-q}.
\end{eqnarray*}
Observing that $B_M\geq C_M$, and that $p-q\frac{mp} 2 = q$, we inject the above inequality in the estimate \iref{estimabove}, which yields
$$
e_{T}(f)^p \leq s^{mp} \left(K_M(\pi_z) + B_M \omega(r)\right)^q |T|.
$$
Averaging on $z\in T$, we obtain
$$
e_{T}(f)^p \leq s^{mp} \int_T  \left(K_M(\pi_z) + B_M \omega(r)\right)^q dz.
$$
Adding up contributions from all triangles in $\cT_s$, we find
$$
 e_{\cT_s}(f)^p =  \sum_{T\in \cTr} e_{T}(f)^p + \sum_{T\in \cTb} e_{T}(f)^p \leq s^{mp} \int_\Omega \left(K_M(\pi_z) + B_M \omega(r)\right)^q dz + C^*_{\rm bd} s^{mp+1}
$$
Combining this with the estimate \iref{cardcTsUpper} we obtain,
$$
e_{\cT_s} \#(\cT_s)^{m/2} \leq \left(\int_\Omega \left(K_M(\pi_z) + B_M \omega(r)\right)^q dz+C^*_{bd} s\right)^{\frac 1 p} 
\left(\int_\Omega \left(K_M(\pi_z)+3B_M \omega(r)\right)^q dz +C_{\rm bd} s \right)^{\frac m 2}
$$
and therefore, since $\frac 1 q = \frac m 2+\frac 1 p$,
$$
\limsup_{s\to 0}\(\#(\cT_s)^{m/2} e_{\cT_s}\) \leq \left(\int_\Omega \left(K_M(\pi_z)+3B_M \omega(r)\right)^q dz \right)^{\frac 1 q}.
$$
It is now time to observe that for fixed $M$,
$$
\lim_{r\ra 0} \int_\Omega \left(K_M(\pi_z) + 3B_M \omega(r)\right)^q dz = \int_\Omega K_M^q(\pi_z) dz,
$$
and that 
$$
\lim_{M\ra +\infty} \int_\Omega K_M^q(\pi_z)dz = \int_\Omega K^q(\pi_z) dz.
$$
Therefore, for all $\e>0$, we can choose $M$ sufficiently large and $r$
sufficiently small, such that
$$
\limsup_{s\to 0}\(\#(\cT_s)^{m/2} e_{\cT_s}\)\leq \(\int_\Omega K^q(\pi_z) dz\)^{\frac 1 q}+\e.
$$
This gives us the announced statement of Theorem \ref{optitheorem}, by defining
$$
s_N:=\min\{s>0 \sep \#(T_s)\leq N\},
$$
and by setting $\cT_N=\cT_{s_N}$.
\nl
\nl
The adaptation of the above proof in the case $p=\infty$ is not straightforward
due to the fact that the contribution to the error of $\cTb$ is not anymore neglectible
with respect to the contribution of $\cTr$. For this reason, one needs to modify
the construction of $\cTb$. Here, we provide a simple construction
but for which the resulting triangulation $\cT_s$ is non-conforming, as we do not know how to produce a satisfying conforming triangulation.

More precisely, we define $\cTr$ in a similar way as for $p<\infty$,
and add to the construction of $\cTb$ a post processing step
in which each triangle is splitted in $4^j$ similar triangles according
to the midpoint rule. Here we take for $j$ the smallest integer
which is larger than $-\frac {\log s}{4\log 2}$. With such an additional splitting,
we thus have
$$
\max_{T\in \cTb}{\rm diam}(T)\leq s^{\frac 1 4}\max_{R\in\cR}{\diam}(sT_R)
=C_as^{1+\frac 1 4}.
$$
The contribution of $\cTb$ to the $L^\infty$ interpolation error is bounded by 
$$
e_{\cTb}(f) \leq C_0C_1\max_{T\in \cTb}{\rm diam}(T)^m \leq C^*_{\rm bd}s^{\frac {5m} 4},
$$
with $C^*_{\rm bd}:=C_0C_1C_a^m$. We also have
$$
\#(\cTb)\leq C_{\rm bd} s^{-3/2},
$$
which remains neglectible compared to $s^{-2}$. We therefore obtain
\be
\label{cardcTsInf}
\#(\cT_s)\leq s^{-2}\left(\int_\Omega (K_M(\pi_z)+3B_M \omega(r))^{\frac 2 m} dz +C_{\rm bd} s^{1/2} \right)
\ee
Moreover, if $T\in \cTr$ and $T\subset R\in \cR$, we have according to the estimate \iref{localabove}
$$
e_{T}(f) \leq \left(K_M(\pi_{b_R}) + B_M \omega(\max\{r,C_A s\})\right) |T|^{\frac{m} 2}.
$$
By construction $|T| = s^2 (K_M(\pi_{b_R})+2B_M \omega(r))^{-2/m}$. This implies $e_{T}(f) \leq s^m$ when $C_A s\leq r$. Therefore
$$
e_{\cT_s}(f) = \max\{e_\cTr,e_\cTb\} \leq s^m\max\{1, C^*_{\rm bd} s^{\frac m 4}\}.
$$
Combining this estimate with \iref{cardcTsInf} yields
$$
\limsup_{s\to 0}\(\#(\cT_s)^{m/2} e_{\cT_s}\) \leq \left(\int_\Omega \left(K_M(\pi_z)+3B_M \omega(r)\right)^{\frac  2 m} dz \right)^{\frac  m 2},
$$
and we conclude the proof in a similar way as for $p<\infty$.

\section{The shape function and the optimal metric for linear and quadratic elements}

This section is devoted to linear ($m=2$) and quadratic ($m=3$) elements,
which are the most commonly used in practice.
In these two cases, we are able to derive an exact expression for $K_{m,p}(\pi)$ in terms
of the coefficients of $\pi$. Our analysis also gives us access to the distorted metric
which characterizes the optimal mesh.
While the results concerning linear elements
have strong similarities with those of \cite{B}, those
concerning quadratic elements are to our knowledge
the first of this kind, although \cite{C1} analyzes a similar setting. 

\subsection{Exact expression of the shape function}

In order to give the exact expression of $K_{m,p}$, we define the determinant
of an homogeneous quadratic polynomial by
$$
\det (ax^2+2bxy+cy^2)=ac-b^2,
$$
and the discriminant of an homogeneous cubic polynomial by
$$
\disc(a x^3+ b x^2 y+ c x y^2+ d y^3) = b^2 c^2 - 4 a c^3 - 4 b^3 d + 18 a b c d - 27 a^2 d^2.
$$
The functions $\det$ on $\H_2$ and $\disc$ on $\H_3$ are homogeneous in the sense that
\be
\label{homogDisc}
\det(\lambda \pi) = \lambda^2\det \pi,\quad \disc(\lambda \pi) = \lambda^4  \disc \pi.
\ee
Moreover, it is well known that they obey an invariance property with respect to linear changes of coordinates $\phi$:
\be
\label{invDisc}
\det(\pi\circ \phi) = (\det \phi)^2 \det \pi,\quad \disc( \pi\circ \phi) =  (\det \phi)^6 \disc \pi.
\ee
Our main result relates $K_{m,p}$ to these quantities.
\begin{theorem}
\label{equal23}
We have for all $\pi\in\H_2$,
$$
K_{2,p}(\pi) = \sigma_p(\det\pi) \sqrt{|\det \pi|},
$$
and for all $\pi\in \H_3$,
$$
K_{3,p}(\pi) = \sigma^*_p(\disc \pi) \sqrt[4]{|\disc \pi|},
$$
where $\sigma_p(t)$ and $\sigma_p^*(t)$ are constants that only depend on the sign of $t$.
\end{theorem}

The proof of Theorem \ref{equal23} relies on the possibility of
mapping and arbitrary polynomial $\pi\in \H_2$ such that 
${\rm det}(\pi)\neq 0$ or $\pi\in\H_3$ such that $\disc(\pi)\neq 0$ onto 
two fixed polynomials $\pi_-$ or $\pi_+$ by a linear change of variable
and a sign change.

In the case of $\H_2$, it is well known
that we can choose $\pi_-=x^2-y^2$ and $\pi_+=x^2+y^2$.
More precisely, to  all $\pi\in H_2$, we associate a
symmetric matrix $Q_\pi$ such that $\pi(z)=\<Q_\pi z,z\>$. This
matrix can be diagonalized according to 
$$
Q_\pi=U^\trans 
\left(
\begin{array}{cc}
\lambda_1 & 0\\
0 & \lambda_2
\end{array}
\right)
U,\quad U\in \cO_2, \; \lambda_1,\lambda_2\in\R.
$$
Then, defining the linear transform
$$
\phi_\pi:=U^\trans 
\left(
\begin{array}{cc}
|\lambda_1|^{-\frac 1 2} & 0\\
0 & |\lambda_2|^{-\frac 1 2}
\end{array}
\right)
$$
and $\lambda_\pi={\rm sign}(\lambda_1)\in \{-1,1\}$, 
it is readily seen that 
$$
\lambda_\pi \pi\circ \phi_\pi =
\left\{
\begin{array}{cl}
x^2+y^2 &\text{ if } \det \pi>0\\
x^2-y^2 &\text{ if } \det \pi<0.
\end{array}
\right.
$$

In the case of $\H_3$, a similar result holds, as shown by the following lemma.

\begin{lemma}
\label{lemmaChgVar}
Let $\pi \in \H_3$. There exists a linear transform
$\phi_\pi$ such that
\be
\label{eqChgVar}
\pi\circ \phi_\pi =
\left\{
\begin{array}{cl}
x(x^2-3y^2) &\text{ if } \disc \pi>0\\
x(x^2+3y^2) &\text{ if } \disc \pi<0.
\end{array}
\right.
\ee
\end{lemma}

\proof
Let us first assume that $\pi$ is not divisible by $y$ so that it can be factorized
as 
$$
\pi = \lambda (x-r_1y) (x-r_2 y) (x-r_3 y),
$$
with $\lambda\in\R$ and $r_i\in\C$. If $\disc \pi>0$, then the $r_i$ are real
and we may assume $r_1<r_2<r_3$. Then, defining
$$
\phi_\pi = \lambda (2\disc \pi)^{-1/3} 
\left(\begin{array}{cc} 
r_1 (r_2 + r_3) -2 r_2 r_3 & (r_2-r_3) r_1 \sqrt 3\\
2 r_1-(r_2+r_3) & (r_2-r_3) \sqrt 3,
\end{array}\right).
$$
an elementary computation shows that  $\pi\circ \phi_\pi = x(x^2 - 3 y^2)$.
If $\disc \pi < 0$, then we may assume that $r_1$ is real, $r_2$ and $r_3$
are complex conjugates with ${\rm Im}(r_2)>0$. Then, defining
$$
\phi_\pi = \lambda (2\disc \pi)^{-1/3} 
\left(\begin{array}{cc} 
r_1 (r_2 + r_3) -2 r_2 r_3 & \mi(r_2-r_3) r_1 \sqrt 3\\
2 r_1-(r_2+r_3) & \mi(r_2-r_3) \sqrt 3 
\end{array}\right).
$$
an elementary computation shows that  $\pi\circ \phi_\pi = x(x^2 +3 y^2)$.
Moreover it is easily checked that $\phi_\pi$ has real entries
and is therefore a change of variable in $\R^2$.

In the case where $\pi$ is divisible by $y$, there exists a rotation $U\in \cO_2$
such that $\t \pi:=\pi\circ U$ is not divisible by $y$. By the invariance property
\iref{invDisc} we know that $\disc \pi=\disc \t \pi$. Thus, we reach the same conclusion
with the choice $\phi_\pi :=U\circ \phi_{\t \pi}$.
\sq

\noindent
{\bf Proof of Theorem \ref{equal23}:} for all $\pi\in \H_2$ such that $\det\pi\neq 0$
and for all change of variable $\phi$ and $\lambda\neq 0$,
we may combine the properties of the determinant in \iref{homogDisc} and \iref{invDisc} 
with those of the shape function established in
Proposition \ref{propinvarK}. This gives us
$$
\frac {K_{2,p}(\pi)} {\sqrt{|\det\pi|}}
=\frac {K_{2,p}(\lambda \pi\circ \phi)} {\sqrt{|\det(\lambda \pi\circ\phi)|}}.
$$
Applying this with $\phi=\phi_\pi$ and $\lambda=\lambda_\pi$, we therefore obtain
$$
K_{2,p}(\pi) =
 \sqrt{|\det \pi|} 
\left\{\begin{array}{cc}
K_{2,p}(x^2+y^2) & \text{ if } \det \pi>0,\\
K_{2,p}(x^2-y^2) & \text{ if } \det \pi<0.
\end{array}\right.
$$
This gives the desired result with $\sigma_p(t)=K_{2,p}(x^2+y^2)$ for $t>0$
and $\sigma_p(t)=K_{2,p}(x^2-y^2)$ for $t<0$. In the case
where $\det\pi=0$, then $\pi$ is of the
form $\pi(x,y) = \lambda (\alpha x +\beta y)^2$ and we conclude by
Proposition \ref{vanishprop} that $K_{2,p}(\pi)=0$.

For all $\pi\in \H_3$ such that $\disc\pi\neq 0$, a similar reasoning yields
$$
K_{3,p}(\pi) = \sqrt[4]{|\disc \pi|} 108^{-\frac{1} 4}
\left\{\begin{array}{cc}
K_{3,p}(x(x^2-3y^2)) & \text{ if } \disc \pi>0,\\
K_{3,p}(x(x^2+3y^2)) & \text{ if } \disc \pi<0.
\end{array}\right.,
$$
where the constant $108$ comes from the fact that
$\disc(x(x^2-3y^2)) = - \disc(x(x^2-3y^2)) = 108$. This gives the
desired result with $\sigma_p^*(t)=108^{-\frac{1} 4}K_{3,p}(x(x^2-3y^2))$ for $t>0$
and $\sigma_p^*(t)=108^{-\frac{1} 4}K_{3,p}(x(x^2+3y^2))$ for $t<0$.
In the case where $\disc\pi=0$, then $\pi$ is of the
form $\pi(x,y) =  (\alpha x +\beta y)^2(\gamma x+\delta y)$ and we conclude by
Proposition \ref{vanishprop} that $K_{3,p}(\pi)=0$.\sq 

\begin{remark}
We do not know any simple analytical expression for the constants involved in 
$\sigma_p$ and $\sigma^*_p$, but these can be found by numerical optimization. These constants are known for some special values of $p$ in the case $m=2$, see for example \cite{B}.
\end{remark}

\subsection{Optimal metrics}

Practical mesh generation techniques such as in \cite{Shew,Bois,Peyre,Bamg,Inria} 
are based on the data of a Riemannian metric, by which we mean
a field $h$ of symmetric definite positive matrices
$$
x \in \Omega \mapsto h(x) \in S_2^+.
$$
Typically, the mesh generator takes the metric $h$ as an input and hopefully
returns a triangulation $\cT_h$ adapted to it in the sense that all triangles are close
to equilateral of unit side length with respect to this metric. 
Recently, it has been rigorously proved in \cite{Shew2, Bois} that some algorithms produce bidimensional meshes obeying these constraints, under certain conditions. This must be contrasted with algorithms based on heuristics, such as \cite{Bamg} in two dimensions, and \cite{Inria} in three dimensions, which have been available for some time and offer good performance \cite{A} but no theoretical guaranties.

For a given function $f$ to be approximated, the field of metrics given 
as input should be such that the 
local errors are equidistributed and the aspect ratios are optimal for the 
generated triangulation. Assuming that the error is measured in $X=L^p$
and that we are using finite elements of degree $m-1$, we can construct
this metric as follows, provided that some estimate of $\pi_z = \frac{d^mf(z)}{m!}$ is available
all points $z\in \Omega$.
An ellipse $E_z$ such that $|E_z|$ is equal
or close to
\be
\sup_{E\in \cE, E\subset \Lambda_{\pi_z}} |E|
\label{maxellips}
\ee
is computed, where $\Lambda_{\pi_z}$ is defined as in \iref{lambdapi}.
We denote by $h_{\pi_z}\in S_2^+$ the associated symmetric definite
positive matrix such that
$$
E_z=\left\{(x,y) \sep (x,y)^T h_{\pi_z} (x,y) \leq 1\right\}.
$$
Let us notice that the supremum in \iref{maxellips} might not
always be attained or even be finite. This particular case is discussed
in the end of this section. Denoting by $\nu>0$ the desired
order of the $L^p$ error on each triangle, we then define the
metric by rescaling $h_{\pi_z}$ according to
$$
h(z)= \frac 1 {\alpha_z^2} h_{\pi_z}\; \; {\rm where}\;\; 
\alpha_z:=\nu^{\frac {p}{mp+2}} |E_z|^{-\frac {1}{mp+2}}.
$$
With such a rescaling, any triangle $T$ designed by the
mesh generator should be comparable to the ellipse
$z+\alpha_z E_z$ centered around $z$ the barycenter of $T$,
in the sense that
\be
z+c_1\alpha_z E_z \subset T\subset z+c_2\alpha_z E_z,
\ee
for two fixed constants $0<2c_1\leq c_2$ independent of $T$ (recall that for any ellipse $E$
there always exist a triangle $T$ such that $E\subset T \subset 2E$).

Such a triangulation heuristically 
fulfills the desired properties of optimal aspect ratio and 
error equidistribution when the level of refinement
is sufficiently small. Indeed, we then have
\begin{eqnarray*}
e_{m,T}(f)_p & \approx & e_{m,T}(\pi_z)_p \\
&=& \|\pi_z - I_{m,T}\pi_z\|_{L^p(T)},\\
&\sim & |T|^{\frac 1 p}\|\pi_z - I_{m,T}\pi_z\|_{L^\infty(T)},\\
&\sim & |T|^{\frac 1 p} \|\pi_z\|_{L^\infty(T)},\\
&\sim &  |\alpha_z E_z|^{\frac 1 p} \|\pi_z\|_{L^\infty(\alpha_z E_z)},\\
&=&	 \alpha_z^{m+\frac 2 p} |E_z|^{\frac 1 p} \|\pi_z\|_{L^\infty(E_z)},\\
&= & \nu,
\end{eqnarray*}
where we have used the fact that $\pi_z\in\H_m$.

Leaving aside these heuristics on error estimation and mesh generation, 
we focus on the main computational issue in the design of the metric $h(z)$,
namely the solution to the problem \iref{maxellips}: to any given $\pi\in \H_m$,
we want to associate $h_\pi\in S_2^+$ such that
the ellipse $E_\pi$ defined by $h_\pi$ has area equal
or close to $\sup_{E\in \cE, E\subset \Lambda_{\pi}} |E|$.

When $m=2$ the computation of the optimal matrix $h_\pi$ 
can be done by elementary algebraic means. In fact, as it will
be recalled below, $h_\pi$ is simply the absolute value
(in the sense of symmetric matrices) of the symmetric
matrix $[\pi]$ associated to the quadratic form $\pi$. These facts
are well known and used in mesh generation algorithms
for $\P_1$ elements.

When $m\geq 3$ no such algebraic derivation of $h_\pi$ from
$\pi$ has been proposed up to now and current 
approaches instead consist in numerically solving 
the optimization problem \iref{OptimEll}, see \cite{C3}. 
Since these computations have to be done extremely frequently 
in the mesh adaptation process, a simpler algebraic procedure 
is highly valuable. In this section, we propose a simple and 
algebraic method in the case $m=3$, corresponding to quadratic elements. 
For purposes of comparison the results already known in the case $m=2$ are recalled.

\begin{prop}
\label{propEllipseMax}
\begin{enumerate}
\item
Let $\pi\in \H_2$ be such that $\det(\pi)\neq 0$, and consider its associated $2\times 2$ matrix which can be 
written as
$$
[\pi] = U^\trans 
\left(
\begin{array}{cc}
\lambda_1 & 0\\
0 & \lambda_2
\end{array}
\right)
U,\quad U\in \cO_2.
$$
Then, an ellipse of maximal volume inscribed in $\Lambda_\pi$
is defined by the matrix
$$
h_\pi = U^\trans 
\left(
\begin{array}{cc}
|\lambda_1| & 0\\
0 & |\lambda_2|
\end{array}
\right)
U
$$

\item
Let $\pi \in \H_3$ be such that $\disc \pi > 0$, and $\phi_\pi$ a matrix satisfying \iref{eqChgVar}. Define
\be
\label{hPiPos}
h_\pi =(\phi_\pi^{-1})^\trans \phi_\pi^{-1}.
\ee

Then $h_\pi$ defines an ellipse of maximal volume inscribed in $\Lambda_\pi$. Moreover $\det h_\pi =\frac {2^{-2/3}} 3 (\disc \pi)^{\frac 1 3}$. 

\item
Let $\pi\in \H_3$ be such that $\disc \pi < 0$, and $\phi_\pi$ a matrix satisfying \iref{eqChgVar}. Define  
$$
h_\pi = 2^{\frac 1 3} (\phi_\pi^{-1})^\trans \phi_\pi^{-1}.
$$
Then $h_\pi$ defines an ellipse of maximal volume inscribed in $\Lambda_\pi$. Moreover $\det h_\pi = \frac 1 3 |\disc \pi|^{\frac 1 3}$ .
\end{enumerate}
\end{prop}

\begin{figure}
	\centering
		\includegraphics[width=5cm,height=5cm]{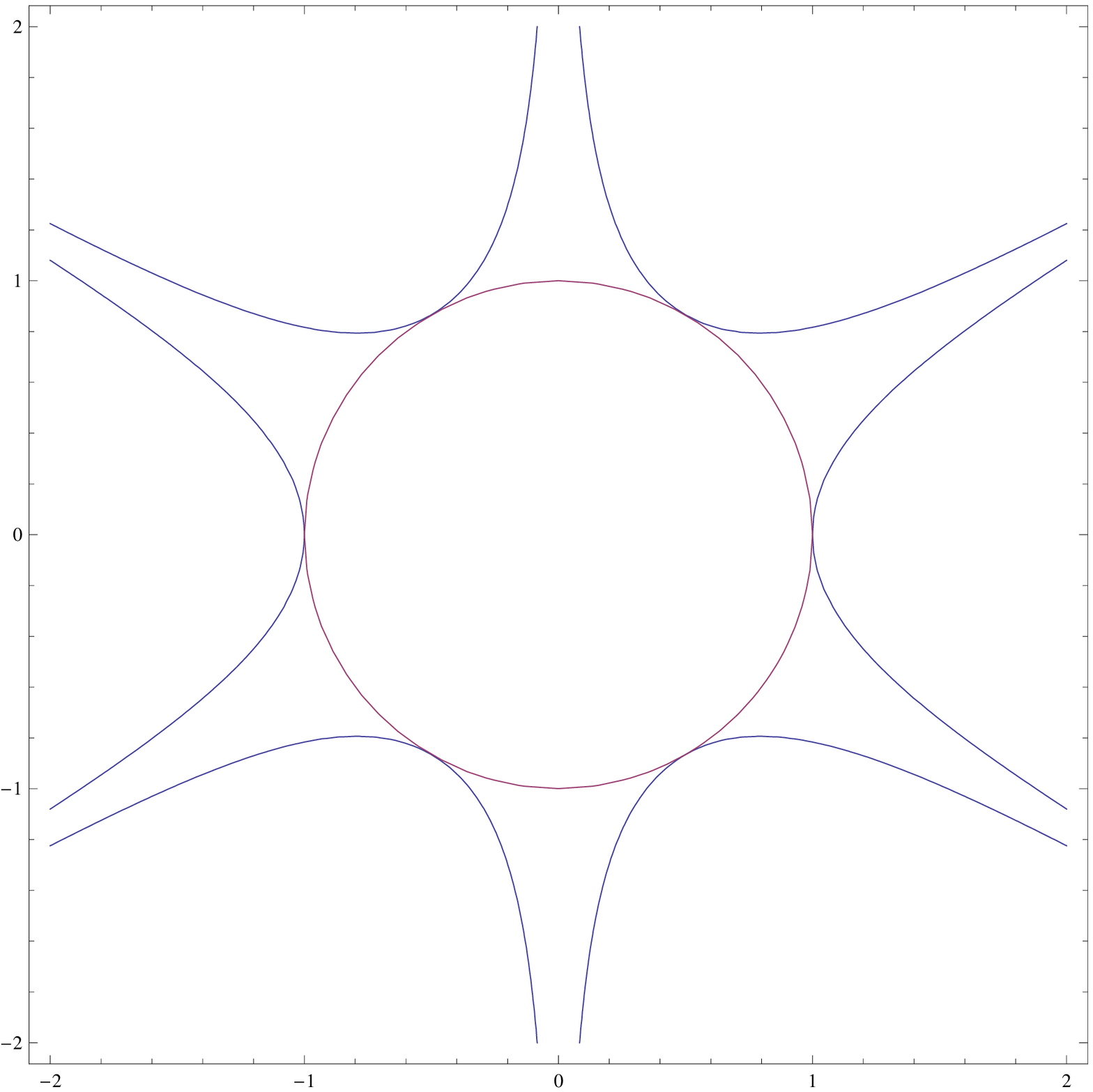}
		\includegraphics[width=5cm,height=5cm]{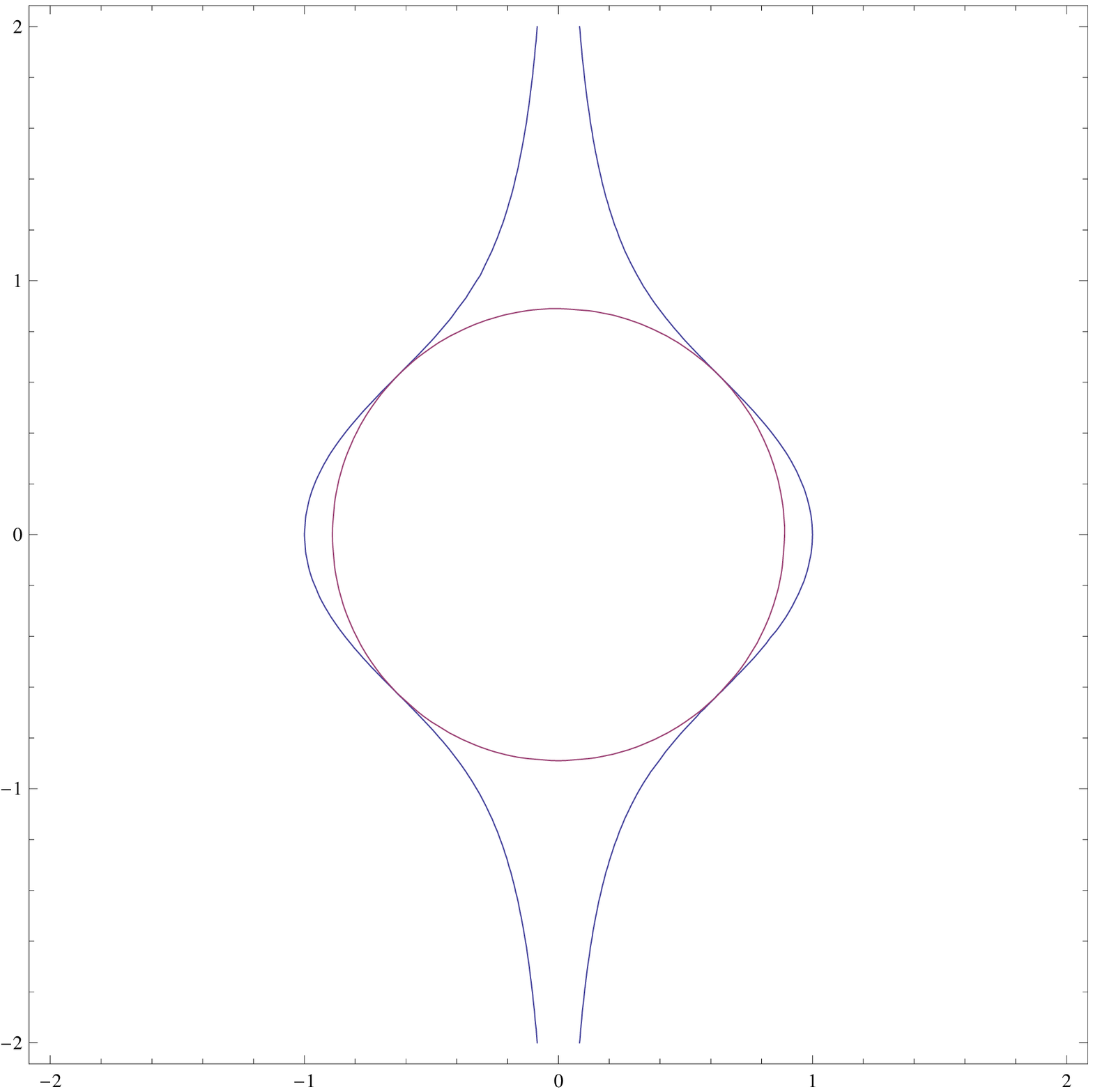}
	\caption{Maximal ellipses inscribed in $\Lambda_\pi$, $\pi = x(x^2-3y^2)$ or $\pi = x(x^2+3y^2)$.}
\end{figure}

\proof
Clearly, if the matrix $h_\pi$ defines an ellipse of maximal volume in the set $\Lambda_\pi$, then for any linear change of coordinates $\phi$, the metric $(\phi^{-1})^\trans h_\pi \phi^{-1}$ defines an ellipse of maximal volume in the set $\Lambda_{\pi\circ \phi}$.
When $\pi\in \H_2$, we know that $\lambda_\pi \pi \circ \phi_\pi = x^2+y^2$ when $\det \pi >0$, and $x^2-y^2$ when $\det \pi <0$, where $|\lambda_\pi|=1$. When $\pi \in \H_3$, 
we know from Lemma \ref{lemmaChgVar} that $\pi \circ \phi_\pi = x(x^2-3y^2)$ 
when $\disc \pi>0$ and $x (x^2+3y^2)$ when $\disc \pi <0$. 
Hence it only remains to prove that when $\pi \in \{x^2+y^2,\ x^2-y^2,\ x(x^2-3y^2)\}$, then $h_\pi=\Id$, which means that the disc of radius $1$ is an ellipse of maximal volume inscribed in $\Lambda_\pi$, while when $\pi = x(x^2+3y^2)$ we have $h_\pi=2^{1/3}\Id$.

The case $\pi = x^2+y^2$ is trivial. We next concentrate on the case $\pi = x(x^2+3y^2)$,
the treatment of the two other cases being very similar.
Let $E$ be an ellipse included in $\Lambda_\pi$, $\pi = x(x^2+3y^2)$. 
Analyzing the variations of the function $\pi(\cos\theta,\sin\theta)$, it is not hard to see that we can rotate $E$ into another ellipse $E'$, also verifying the inclusion $E'\subset \Lambda_\pi$, and which principal axes are $\{x=0\}$ and $\{y=0\}$.
We therefore only need to consider ellipses of the form $k x^2 + h y^2 \leq 1$. 
For a given value of $h$, we denote by $k(h)$ the minimal value of $k$ for which this ellipse is included in $\Lambda_\pi$. Clearly the boundary of the ellipse, defined by $k(h) x^2 + h y^2 = 1$, must be tangent to the curve defined by $\pi(x,y) = 1$ at some point $(x,y)$. This translates into the following system of equations
\be 
\left\{
\begin{array}{ccc}
\pi(x,y) &=& 1,\\
h x^2 + k y^2 &=& 1,\\
k y \partial_x \pi(x,y) - hx \partial_x \pi(x,y) &=& 0.
\end{array}
\right.
\label{tangencyEqn}
\ee
Eliminating the variables $x$ and $y$ from this system, as well as negative
or complex valued solutions, we find that $k(h) = \frac{4+h^3}{3h^2}$ when $h\in (0,2]$, and $k(h) = k(2) = 1$ when $h\geq 2$.
The minimum of the determinant $h k(h) = \frac 1 3 \left(\frac 4 h + h^2\right)$ is attained for $h=2^{\frac 1 3}$. Observing that $k(2^{\frac 1 3}) = 2^{\frac 1 3}$ we obtain as announced $h_\pi = 2^{1/3}\Id$ and that the ellipse of largest area included in $\Lambda_\pi$ is the disc
of equation $2^{1/3}(x^2+y^2)\leq 1$, as illustrated on Figure 2.b.

The same reasoning applies to the other cases. For $\pi = x^2-y^2$ we obtain $k(h) = \frac 1 h$, $h\in (0,\infty)$. In this case the determinant $h k(h)$ is independent of $h$, and we simply choose $h=1 = k(1)$.
For $\pi = x(x^2-3y^2)$ we obtain $k(h) = \frac{4-h^3}{3h^2}$ when $h\in (0,1]$ and $k(h) = k(1) = 1$ when $h>1$. The maximal volume is attained when $h=1$, corresponding to the
unit disc, as illustrated on Figure 2.a.
\sq

\begin{remark}
When $\pi\in \H_3$ and $\disc \pi>0$ a surprising simplification happens : the matrix 
{\rm \iref{hPiPos}} has entries which are symmetric functions of the roots $r_1,r_2,r_3$.
Using the relation between the roots and the coefficients of a polynomial, we find the following expression

$$
\text{If }\pi = a x^3+ 3 b x^2 y+ 3 c x y^2+ d y^3,\text{ then } h_\pi = 2^{-\frac 1 3} 3 (\disc \pi)^{\frac{-1} 3} 
\left(
\begin{array}{cc}
2 (b^2-ac) & bc - ad\\
bc - ad & 2 (c^2-bd)
\end{array}
\right).
$$
This yields a direct expression of the matrix as a function of the coefficients. Unfortunately there is no such expression when $\disc \pi<0$.
\end{remark}

At first sight, Proposition \ref{propEllipseMax} might seem to be a complete solution to the problem of building an appropriate metric for mesh generation. However, some difficulties arise at points $z\in \Omega$ where $\det \pi_z=0$ or $\disc \pi_z=0$.
If $\pi\in \H_2\bs \{0\}$ and $\det \pi = 0$, then up to a linear change of coordinates, and a change of sign, we can assume that $\pi = x^2$. The minimization problem clearly yields the degenerate matrix
$h_\pi = \diag(1,0)$, the $2\times 2$ diagonal matrix with entries $1$ and $0$.
If $\pi\in \H_3\bs \{0\}$ and $\disc \pi = 0$, then up to a linear change of coordinates either $\pi = x^3$ or $\pi = x^2 y$. In the first case the minimization problem gives again $h_\pi = \diag(1,0)$. In the second case a wilder behavior appears, in the sense that minimizing sequences
for the problem \iref{maxellips} are of the type
$h_\pi = \diag(\ve^{-1},\ve^2)$ with $\ve\to 0$. The minimization process 
therefore gives a matrix which is not only degenerate, but also unbounded.

These degenerate cases appear generically, and constitute a problem for
mesh generation since they mean that the adapted triangles are not well defined. 
Current anisotropic mesh generation algorithms for linear elements 
often solve this problem by fixing a small parameter $\delta>0$, and working with the modified matrix
$\tilde h_\pi := h_\pi+\delta\Id$ which cannot degenerate. However this procedure cannot be extended to quadratic elements, since $h_{x^2 y}$ is both degenerate and unbounded.

In the theoretical construction of an optimal mesh which
was discussed in \S 3.2, we  tackled this problem by imposing a bound $M>0$ on the diameter of the triangles. This was the purpose of the modified shape function $K_M(\pi)$ and of the triangle $T_M(\pi)$ of minimal interpolation error among the triangles of diameter smaller than $M$.
We follow a similar idea here, looking for the ellipse of largest area included in $\Lambda_\pi$ with constrained diameter. This provides matrices which are both positive definite and bounded, 
and vary continuously with respect to the data $\pi\in \H_3$. 
The constrained problem, depending on $\alpha>0$, is the following:
\be
\sup \{|E|\sep E\in \cE,\; E\subset \Lambda_\pi \text{ and }\diam E\leq 2\alpha^{-1/2}\},
\label{ellipsconstrained}
\ee
or equivalently
\be
\label{hConstrained}
\inf \{\det H\sep  H\in S_2^+ \;\;{\rm s.t.}\;\; \<Hz,z\> \geq |\pi(z)|^{2/m}, z\in\R^2,\; \text{ and } H\geq \alpha \Id\}.
\ee
We denote by $E_{\pi,\alpha}$ and $h_{\pi,\alpha}$ the solutions
to \iref{ellipsconstrained} and \iref{hConstrained}.
In the remainder of this section, we show that this solution can also be computed 
by a simple algebraic procedure, avoiding any kind of numerical optimization. In the case
where $\pi\in \H_2$, it can easily be checked that
\be
[h_{\pi,\alpha}] = U^\trans 
\left(
\begin{array}{cc}
\max\{|\lambda_1|,\alpha\} & 0\\
0 & \max\{|\lambda_2|,\alpha\},
\end{array}
\right)
U
\label{family2}
\ee
as illustrated on Figure 3. 

\begin{figure}
	\centering
		\includegraphics[width=6cm,height=4cm]{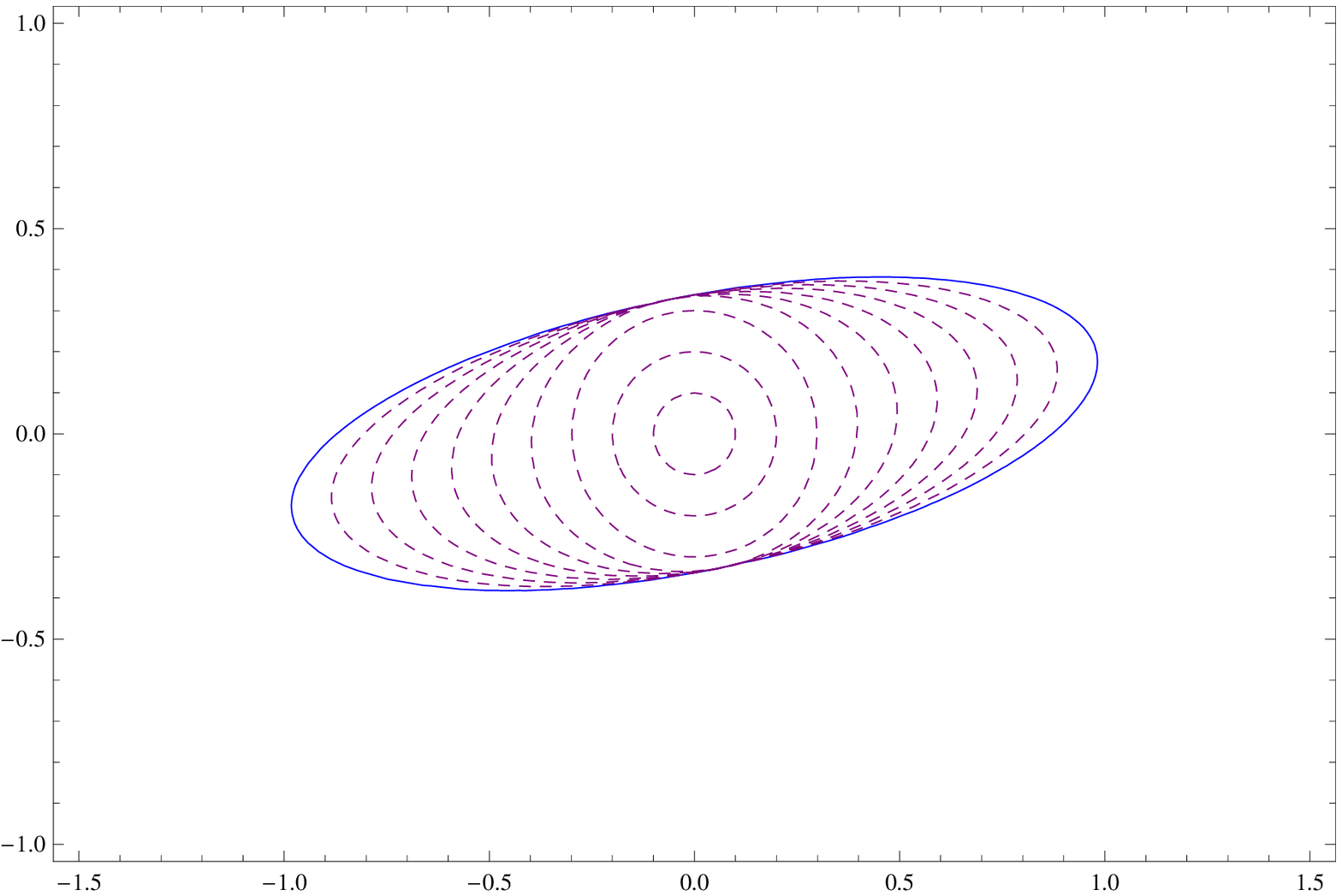}
		\includegraphics[width=6cm,height=4cm]{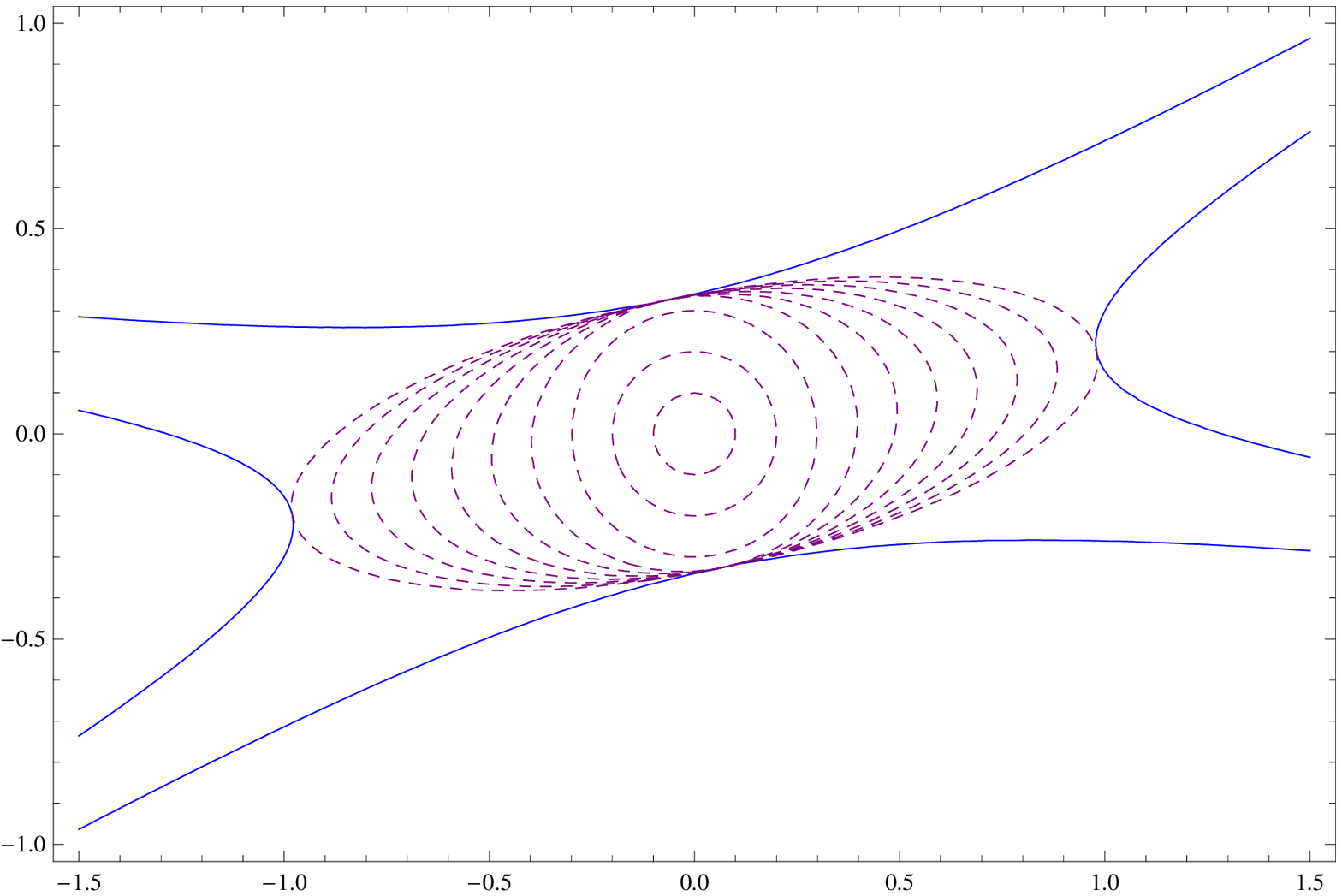}
	\caption{The set $\Lambda_\pi$ (full) and the ellipses $E_{\pi,\alpha}$ (dashed)
    for various values of $\alpha>0$ when $\pi\in \H_2$.}
\end{figure}

When $\pi \in \H_3$, the problem is
more technical, and the matrix $h_{\pi,\alpha}$ takes different forms depending on the value of $\alpha$ and the sign of $\disc \pi$.
In order to describe these different regimes, we introduce three real numbers 
$0\leq \beta_\pi\leq \alpha_\pi\leq \mu_\pi$ and a matrix 
$U_\pi\in \cO_2$ which are defined as follow. We first define $\mu_\pi$ by
$$
\mu_\pi^{-1/2}:=\min\{\|z\|\sep |\pi(z)|=1\},
$$ 
the radius of the largest
disc $D_\pi$ inscribed in $\Lambda_\pi$. For $z_\pi$ such that $|\pi(z_\pi)|=1$
and $\|z_\pi\|=\mu_\pi^{-1/2}$, we define $U_\pi$
as the rotation which maps $z_\pi$ to the vector $(\|z_\pi\|,0)$.
We then define $\alpha_\pi$ by
$$
2\alpha_\pi^{-1/2}:=\max \{{\rm diam}(E)\sep E\in \cE\; ;\; D_\pi\subset E\subset\Lambda_\pi\},
$$
the diameter of the largest ellipse inscribed in $\Lambda_\pi$
and containing the disc $D_\pi$. In the case where $\pi$ is
of the form $(ax+by)^3$, this ellipse is infinitely long
and we set $\alpha_\pi=0$. We finally define $\beta_\pi$ by
$$
2\beta_\pi^{-1/2}:={\rm diam}(E_\pi),
$$ 
where $E_\pi$ is the optimal
ellipse described in Proposition \ref{propEllipseMax}. In the case
where $\disc\pi=0$, the ``optimal ellipse'' is infinitely long
and we set $\beta_\pi=0$. It is readily seen that 
$0\leq \beta_\pi\leq \alpha_\pi\leq \mu_\pi$.

All these quantities can be algebraically computed
from the coefficients of $\pi$ by solving equations of degree
at most $4$, as well as the other quantities involved in the
description of the optimal $h_{\pi,\alpha}$ and $E_{\pi,\alpha}$
in the following result.

\begin{prop}
\label{family3}
For $\pi\in \H_3$ and $\alpha>0$, the matrix $h_{\pi,\alpha}$ and ellipse $E_{\pi,\alpha}$ are described as follows.
\begin{enumerate}
\item
If $\alpha\geq \mu_\pi$, then $h_{\pi,\alpha}=\alpha\Id$
and $E_{\alpha,\pi}$ is the disc of radius $\alpha^{-1/2}$.
\item
If $\alpha_\pi\leq \alpha\leq \mu_\pi$, then 
\be
\label{hBigAlpha}
h_{\pi,\alpha} = U_\pi^\trans 
\left(
\begin{array}{cc}
\mu_\pi & 0\\
0 & \alpha
\end{array}
\right)
U_\pi,
\ee
and $E_\alpha$ is the ellipse of diameter $2\alpha^{-1/2}$
which is inscribed in $\Lambda_\pi$ and contains $D_\pi$.
It is tangent to $\partial\Lambda_\pi$ at the two points $z_\pi$ and $-z_\pi$.
\item
If $\beta_\pi\leq \alpha\leq \alpha_\pi$ then $E_{\pi,\alpha}$ is tangent
to $\partial\Lambda_\pi$ at four points
and has diameter $2\alpha^{-1/2}$. There are at most 
three such ellipses and $E_{\pi,\alpha}$ is the one
of largest area. The matrix $h_{\pi,\alpha}$ has a form 
which depends on the sign of $\disc\pi$.
\nl
(i) If $\disc \pi<0$,
then
$$
h_{\pi,\alpha} = (\phi_\pi^{-1})^\trans 
\left(
\begin{array}{cc}
\lambda_\alpha & 0\\
0 & \frac{4+ \lambda_\alpha^3}{3\lambda_\alpha^2}
\end{array}
\right)
\phi_\pi^{-1}
$$
where $\phi_\pi$ is the matrix defined in Proposition {\rm \ref{propEllipseMax}} 
and $\lambda_\alpha$ determined by $\det(h_{\pi,\alpha}-\alpha\Id)=0$.
\nl
(ii)
If $\disc \pi >0$, then
$$
h_{\pi,\alpha} = (\phi_\pi^{-1})^\trans V^\trans
\left(
\begin{array}{cc}
\lambda_\alpha & 0\\
0 & \frac{4-\lambda_\alpha^3}{3\lambda_\alpha^2}
\end{array}
\right)
V
\phi_\pi^{-1}
$$
where $\phi_\pi$ and $\lambda_\alpha$ are given as in the case $\disc\pi<0$ and
where $V$ is chosen between the three rotations by $0$, $60$ or $120$ degrees
so to maximize $|E_{\alpha,\pi}|$.\nl
(iii)
If $\disc \pi = 0$ and $\alpha_\pi>0$, then there exists a linear change 
of coordinates $\phi$ such that $\pi\circ\phi =x^2 y$ and we have 
$$
h_{\pi,\alpha} = (\phi^{-1})^\trans
\left(
\begin{array}{cc}
\lambda_\alpha & 0\\
0 & \frac 4 {27\lambda_\alpha^2}
\end{array}
\right)
\phi^{-1}
$$
where $\lambda_\alpha$ is determined by $\det(h_{\pi,\alpha}-\alpha\Id)=0$.
\item
If $\alpha\leq \beta_\pi$, then $h_{\pi,\alpha}=h_\pi$ and
$E_{\pi,\alpha}=E_\pi$ is the solution of the unconstrained problem.
\end{enumerate}
\end{prop}

\proof
See Appendix.
\sq

\begin{figure}
	\centering
		\includegraphics[width=4cm,height=4cm]{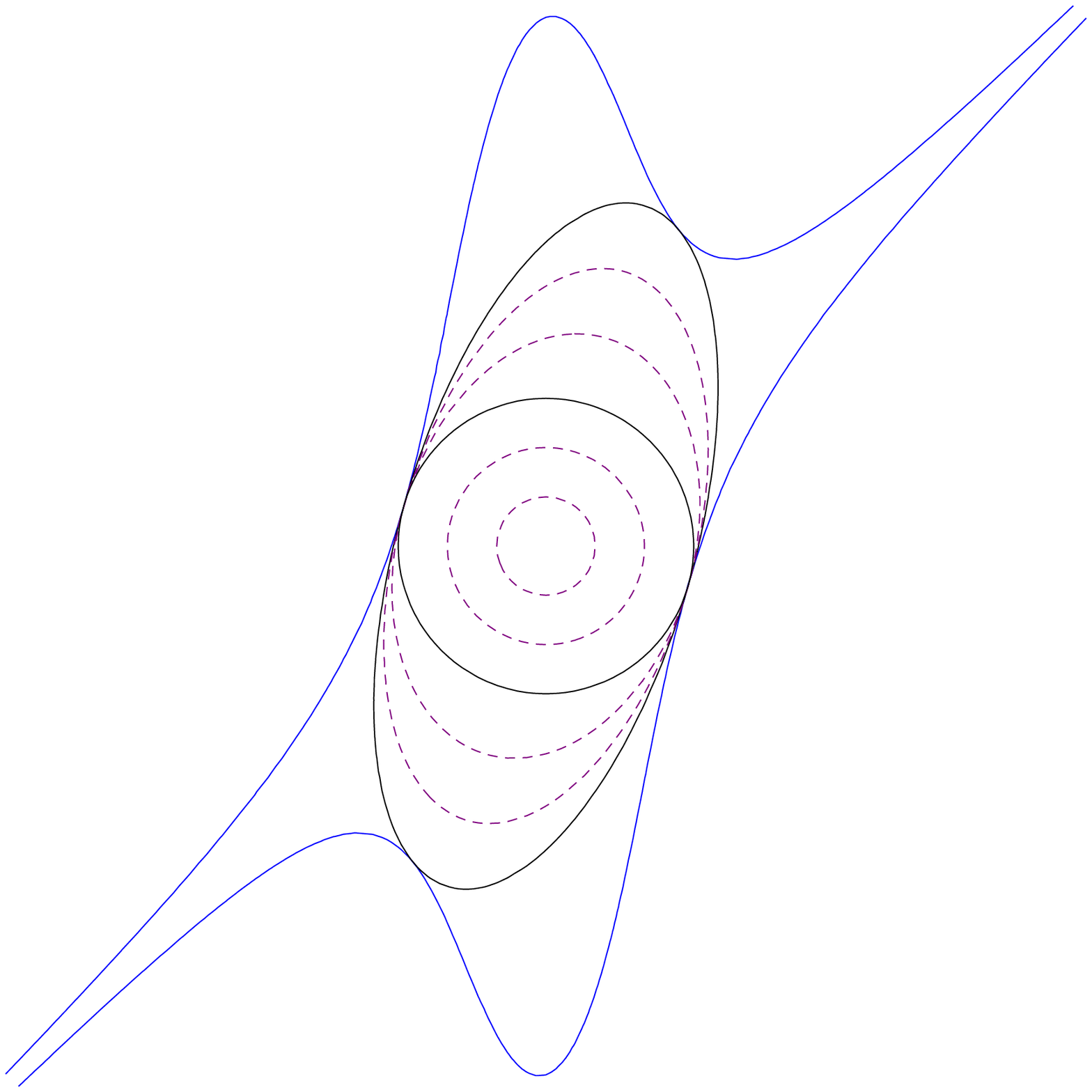}
		\includegraphics[width=4cm,height=4cm]{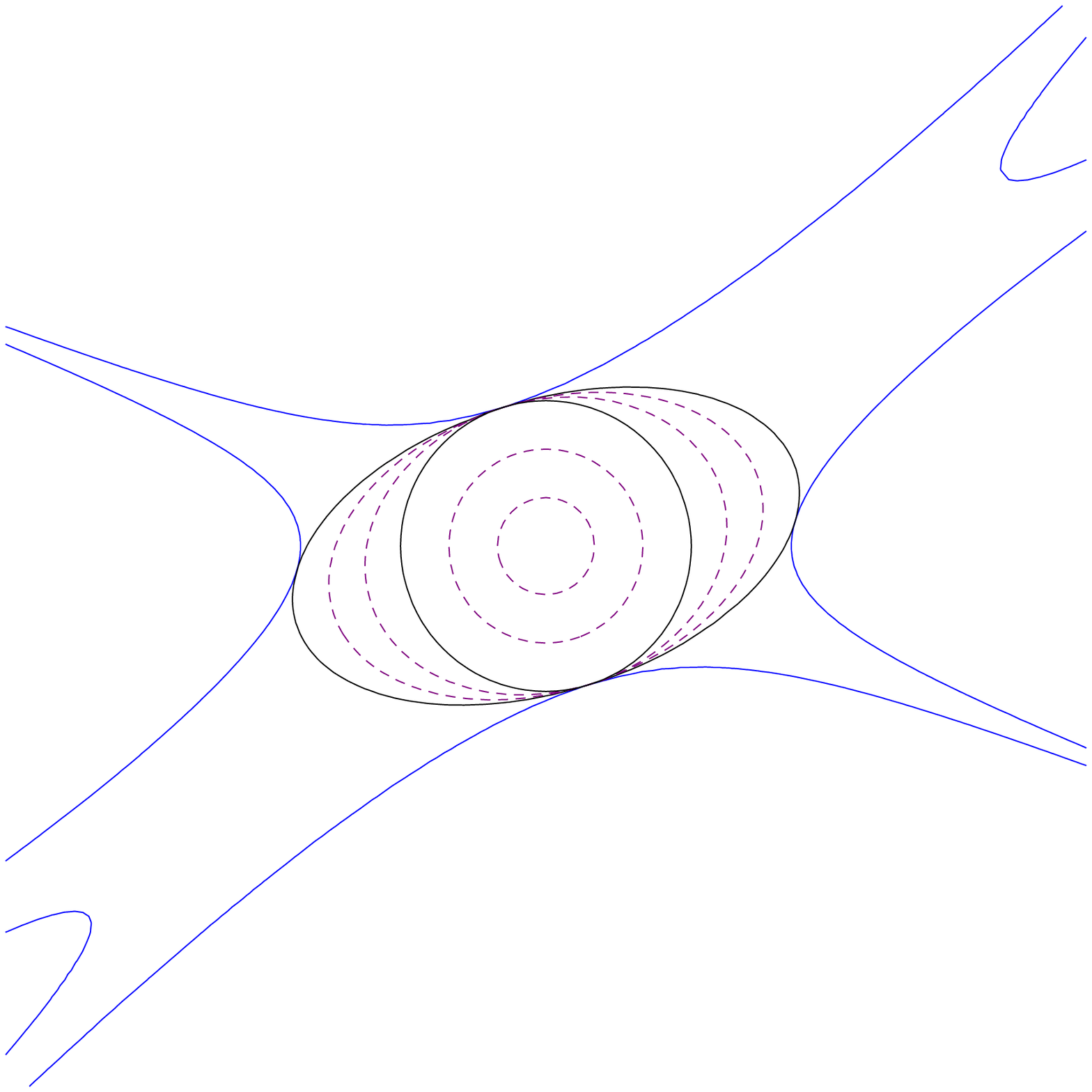}
	\caption{The set $\Lambda_\pi$ (full), the disc $E_{\pi,\mu_\pi}=D_\pi$ (full), 
    the ellipse $E_{\pi,\alpha_\pi}$ (full), and the ellipses $E_{\pi,\alpha}$ (dashed)
    for various values of $\alpha>0$ when $\pi\in \H_3$ and $\alpha \in (\alpha_\pi,\infty)$.
    Left: $\disc\pi <0$. Right: $\disc\pi >0$}
\end{figure}

\begin{figure}
	\centering
		\includegraphics[width=4cm,height=4cm]{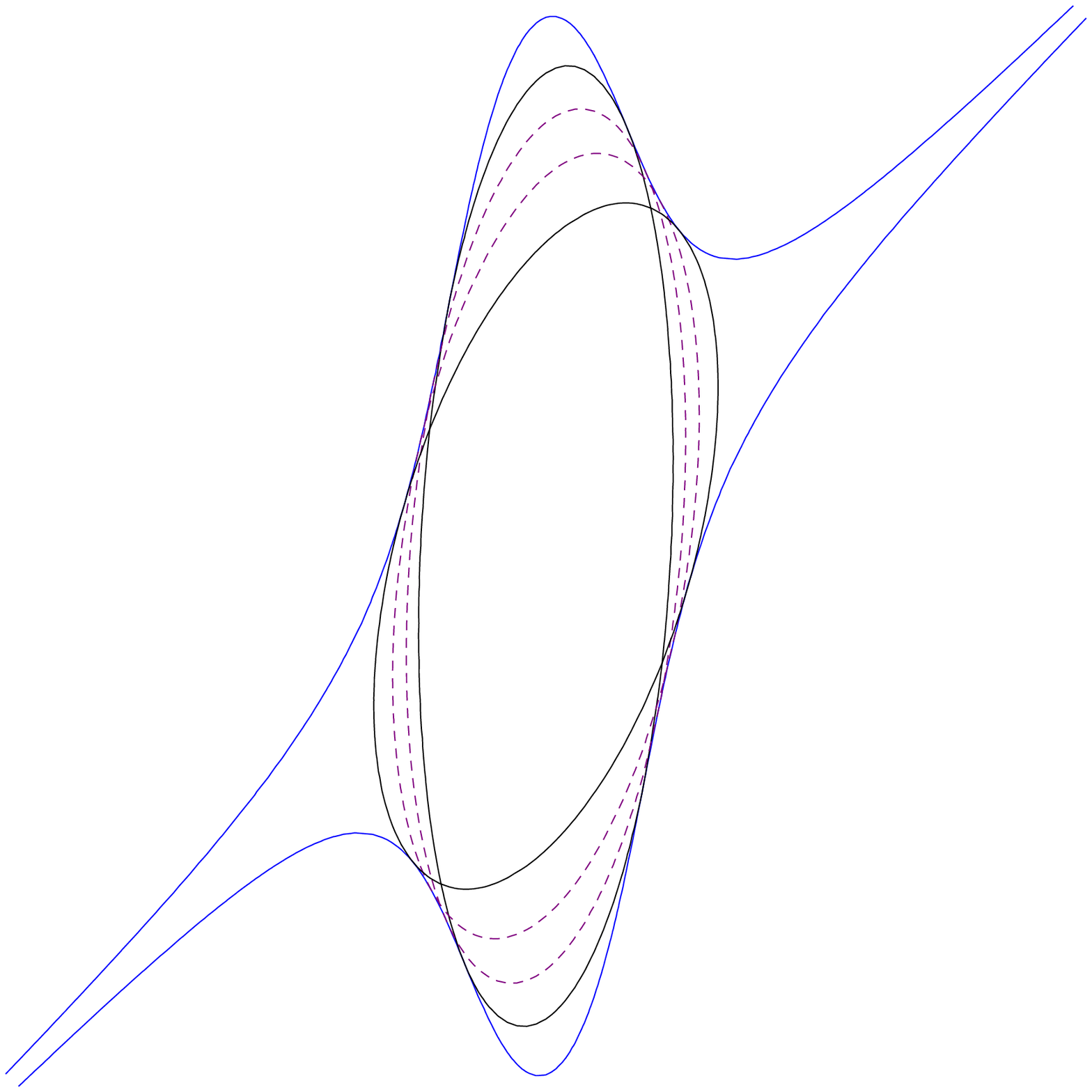}
		\includegraphics[width=4cm,height=4cm]{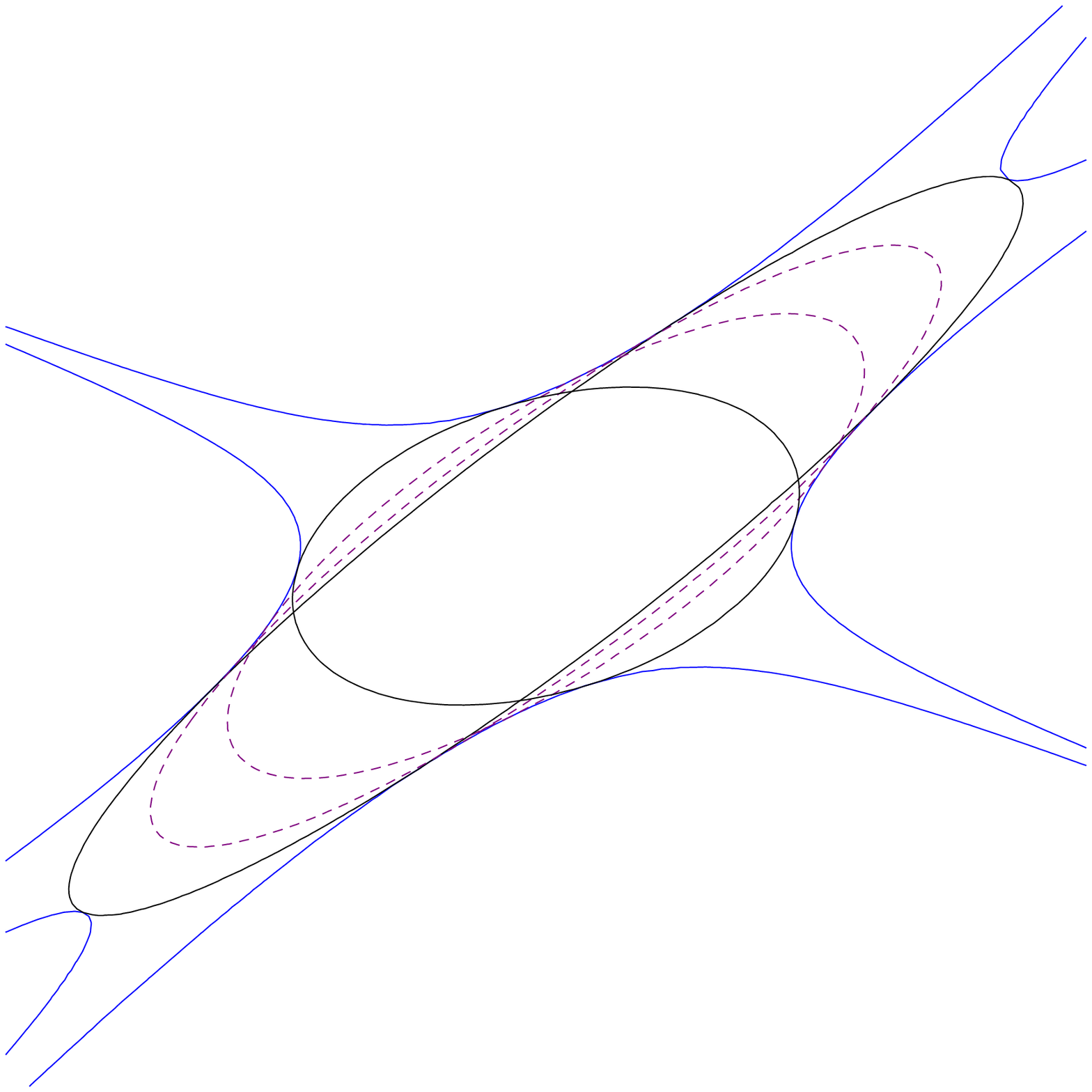}
	\caption{The set $\Lambda_\pi$ (full), the ellipse $E_{\pi,\alpha_\pi}$ (full), 
    the ellipse $E_{\pi,\beta_\pi}=E_\pi$ (full), and the ellipses $E_{\pi,\alpha}$ (dashed)
    for various values of $\alpha>0$ when $\pi\in \H_3$ and $\alpha \in (\beta_\pi,\alpha_\pi)$.
    Left: $\disc\pi <0$. Right: $\disc\pi >0$}
\end{figure}

\noindent
Figure 4 illustrates the ellipses $E_{\pi,\alpha}$, $\alpha \in (\alpha_\pi,\infty)$ when $\disc \pi>0$ (4.a) or $\disc\pi<0$ (4.b). Figure 5 illustrates the ellipses $E_{\pi,\alpha}$, $\alpha \in (\beta_\pi,\alpha_\pi)$ when $\disc \pi>0$ (5.a) or $\disc\pi<0$ (5.b). Note that when $\alpha\geq \alpha_\pi$,
the principal axes of $E_{\pi,\alpha}$ are independent of $\alpha$ since $U_\pi$
is a rotation that only depends on $\pi$, while these axes generally
vary when $\beta_\pi\leq \alpha \leq \alpha_\pi$, since the matrix $\phi_\pi$ is not
a rotation.

\begin{remark}
\label{overfitting}
For interpolation by cubic or higher degree polynomials ($m\geq 4$), an additional difficulty arises that can be summarized as follows: one should be careful not to ``overfit'' the polynomial $\pi$ with the matrix $h_\pi$. An approach based on exactly solving the optimization problem \iref{maxellips} might indeed
lead to a metric $h(z)$ with unjustified strong variations with respect to $z$ 
and/or bad conditioning, and jeopardize the mesh generation process. 
As an example, consider the one parameter family of polynomials 
$$
\pi_t=x^2y^2+ty^4\in\H_4,\;\; t\in [-1,1].
$$ 
It can be checked that when $t>0$, the supremum
$S_+=\sup_{E\in\cE,E\subset\Lambda_{\pi_t}} |E|$ is finite and independent of $t$, but not attained, and that
any sequence $E_n\subset\Lambda_{\pi_t}$ of ellipses such that $\lim_{n\to\infty} |E_n| = S_+$ becomes infinitely elongated in the $x$ direction, as $n\to\infty$. For $t<0$, the supremum
$S_-=\sup_{E\in\cE,E\subset\Lambda_{\pi_t}} |E|$ is independent of $t$ 
and attained for the optimal ellipse of equation 
$|t|^{-1/2}\frac {\sqrt 2-1} 2 x^2+|t|^{1/2}y^2\leq 1$. This ellipse becomes infinitely elongated
in the $y$ direction as $t\to 0$. This example shows the instability of the optimal matrix $h_\pi$ with respect
to small perturbations of $\pi$. However, for all values of $t\in [-1,1]$, 
these extremely elongated ellipses could be discarded in favor, for example,
of the unit disc $D=\{x^2+y^2\leq 1\}$ which obviously satisfies $D\subset \Lambda_{\pi_t}$ and is a near-optimal
choice in the sense that $2|D|= S_+\leq S_-= |D|\sqrt{2(\sqrt 2+1)}$. 
\end{remark}

\section{Polynomial equivalents of the shape function in higher degree}

In degrees $m\geq 4$, we could not find analytical expressions of $K_{m,p}$ or $K_m^\cE$, and do not expect them to exist. However, equivalent quantities with analytical expressions
are available, under the same general form as in Theorem \ref{equal23}: the root of a polynomial in the coefficients of the polynomial $\pi\in\H_m$. This result improves on the analysis of \cite{C2}, where a similar setting is studied.

In the following, we say that a function $\rm \b R$ is a polynomial on $\H_m$ if there exists a polynomial $P$ of $m+1$ variables such that for all $(a_0,\cdots, a_m)\in \R^{m+1}$,
$$
{\rm \b R}\left(\sum_{i=0}^m a_i x^i y^{m-i}\right) := P(a_0,\cdots, a_m) ,
$$
and we define $\deg {\rm \b R} := \deg P$.

The object of this section is to prove the following theorem 
\begin{theorem}
\label{thequiv}
For all degree $m\geq 2$, there exists a polynomial $\Kpol_m$ on $\H_m$, 
and a constant $C_m>0$ such that for all $ \pi \in \H_m$, and all $1\leq p \leq \infty$
$$
\frac 1 {C_m} \sqrt[r_m]{\Kpol_m(\pi)} \leq K_{m,p}(\pi )\leq C_m \sqrt[r_m]{\Kpol_m(\pi)},
$$
where $r_m = \deg \Kpol_m$.
\end{theorem}

Since for fixed $m$ all functions $K_{m,p}$, $1\leq p\leq \infty$, are equivalent on $\H_m$, there is no need to keep track of the exponent $p$ in this section and we use below the notation $K_m = K_{m,\infty}$. In this section, please do not confuse the functions $K_m$ and $\Kpol_m$, as well as the polynomials $Q_d$ and $\b Q_d$ below, which notations are only distinguished by their case.

Theorem \ref{thequiv} is a generalization of Theorem \ref{equal23}, and the polynomial $\Kpol_m$ involved should be seen as a generalization of the determinant on $\H_2$, and of the discriminant on $\H_3$.
Let us immediately stress that the polynomial $\Kpol_m$ is not unique. In particular, we shall propose
two constructions that lead to different $\Kpol_m$ with different degree $r_m$.
Our first construction is simple and intuitive, but leads to a polynomial 
of degree $r_m$ that grows quickly with $m$. Our second construction
uses the tools of Invariant Theory to provide a polynomial
of much smaller degree, which might be more useful in practice.

We first recall that there is a strong connection between the roots of a polynomial in $\H_2$ or $\H_3$ and its determinant or discriminant.
\begin{eqnarray*}
\det\left(\lambda \prod_{1\leq i\leq 2} (x-r_i y)\right) & = & \frac {-1} 4 \lambda^2(r_1-r_2)^2,\\ 
\disc\left(\lambda\prod_{1\leq i\leq 3} (x-r_i y)\right)  & = & \lambda^4 (r_1-r_2)^2(r_2-r_3)^2(r_3-r_1)^2.
\end{eqnarray*}

We now fix an integer $m>3$. Observing that these expressions are a ``cyclic'' product of the squares of differences of roots, we define
$$
\cyc(\lambda,r_1,\cdots,r_m) := \lambda^4 (r_1-r_2)^2\cdots (r_{m-1}-r_m)^2 (r_m-r_1)^2.
$$
Since $m>3$, this quantity is not invariant anymore under reordering of the $r_i$. For any
positive integer $d$, we
introduce the symmetrized version of the $d$-powers of the cyclic product
$$
Q_d(\lambda,r_1,\cdots,r_m) := \sum_{\sigma\in \Sigma_m} \cyc(\lambda,r_{\sigma_1},\cdots,
r_{\sigma_m})^d,
$$
where $\Sigma_m$ is the set of all permutations of $\{1,\cdots, m\}$. 

\begin{prop}
For all $d>0$ there exists a homogeneous polynomial $\b Q_d$ of degree $4 d$ on $\H_m$, with integer coefficients, and such that
$$
\text{If } \pi = \lambda \prod_{i=1}^m (x-r_i y)
\text{ then } \b Q_d \left(\pi \right) = Q_d(\lambda,r_1,\cdots,r_m).
$$ 
In addition, $\b Q_d$ obeys the invariance property
\be
\label{invprop}
\b Q_d(\pi\circ\phi) = (\det \phi)^{2md} \b Q_d(\pi).
\ee
\end{prop}

\proof
We denote by $\sigma_i$ the elementary symmetric functions in the $r_i$, in such way that
$$
\prod_{i=1}^m (x-r_i y) = x^m - \sigma_1 x^{m-1} y+ \sigma_2 x^{m-2} y^2 -\cdots +(-1)^m \sigma_m y^m.
$$
A well known theorem of algebra (see e.g. chapter IV.6 in \cite{Lang}) asserts that any symmetrical polynomial in the $r_i$, can be reformulated as a polynomial in the $\sigma_i$. Hence for any $d$ there exists a polynomial $\tilde Q_d$ such that
$$
Q_d(1,r_1,\cdots,r_m) = \tilde Q_d(\sigma_1,\cdots,\sigma_m).
$$
In addition it is known that the total degree of $\tilde Q_d$ is the partial degree of $Q_d$ in the variable $r_1$, in our case $4d$, and that $\tilde Q_d$ has integer coefficients since $Q_d$ has.

Given a polynomial $\pi\in H_m$ not divisible by $y$, we write it under the two equivalent forms
$$
\pi = a_0 x^m + a_1 x^{m-1}y+\cdots + a_m y^m = \lambda \prod_{i=1}^m(x-r_i y).
$$
clearly $a_0 = \lambda$ and $\sigma_i = (-1)^i \frac{a_i}{a_0}$. It follows that 
$$
Q_d(\lambda,r_1,\cdots,r_m)= \lambda^{4d}\tilde Q_d(\sigma_1, \cdots,\sigma_m) = a_0^{4d} \tilde Q_d\left(\frac{-a_1}{a_0},\cdots,\frac{(-1)^m a_m}{a_0}\right)
$$
Since $\deg \tilde Q_d = 4d$, the negative powers of $a_0$ due to the denominators are cleared 
by the factor $a_0^{4d}$ and the right hand side is thus 
a polynomial in the coefficients $a_0,\cdots,a_m$ that we denote by $\b Q_d(\pi)$.

We now prove the invariance of $\b Q_d$ with respect to linear changes of coordinates, this proof is adapted from \cite{Hilbert}. 
By continuity of $\b Q_d$, it suffices to prove this invariance property for pairs $(\pi,\phi)$ such that $\phi$ is an invertible linear change of coordinates, and neither $\pi$ or $\pi\circ \phi^{-1}$ is divisible by $y$.

Under this assumption, we observe that if $\pi = \lambda \prod_{i=1}^m(x-r_i y)$ and $\phi = \left(\begin{array}{cc} \alpha &\beta\\ \gamma &\delta \end{array}\right)$, then $\pi \circ \phi^{-1} = \tilde \lambda \prod_{i=1}^m (x-\tilde r_i y)$ where

$$
\tilde \lambda = \lambda (\det\phi)^{-m} \prod_{i=1}^m(\gamma+\delta r_i) \text{ and } \tilde r_i =\frac {\alpha r_i+\beta }{\gamma r_i+\delta}.
$$
Observing that 
$$
\tilde r_i - \tilde r_j = \frac{\det\phi}{(\gamma r_i+\delta)(\gamma r_j+\delta)} (r_i-r_j),
$$
it follows that 
$$
\cyc(\tilde \lambda,\tilde r_1,\cdots,\tilde r_m) = (\det \phi)^{-2m} \cyc(\lambda,r_1,\cdots,r_m).
$$
The invariance property \iref{invprop} follows readily.
\sq
\nl
We now define $r_m = 2 {\rm lcm}\{\deg \b Q_d\sep  1\leq d \leq m!\}$
where ${\rm lcm}\{a_1,\cdots,a_k\}$ stands for the lowest common multiple of 
$\{a_1,\cdots,a_k\}$, and we consider the following polynomial on $\H_m$: 

$$
\Kpol_m := \sum_{d=1}^{m!} \b Q_d^{\frac{r_m}{\deg \b Q_d} }.
$$
Clearly $\Kpol_m$ has degree $r_m$ and obeys the invariance property $\Kpol_m(\pi\circ\phi) = (\det\phi)^{\frac{r_m m} 2} \Kpol_m(\pi)$. 
\begin{lemma}
Let $\pi\in\H_m$. If $\Kpol_m(\pi) = 0$ then $K_m(\pi)=0$.
\end{lemma}

\proof
We assume that $K_m(\pi)\neq 0$ and intend to prove that $\Kpol_m(\pi) \neq 0$. Without loss of generality, we may assume that $y$ does not divide $\pi$, since $K_m(\pi\circ U)= K_m(\pi)$ and $\Kpol_m(\pi\circ U)= \Kpol_m(\pi)$ for any rotation $U$. We thus write $\pi = \lambda\prod_{i=1}^m(x-r_i y)$, where $r_i\in\C$. Since $K_m(\pi)\neq 0$, we know
from Proposition \ref{vanishprop} that there is no group of $\mhalf:= \lfloor \frac m 2 \rfloor +1$ equal roots $r_i$.
 
We now define a permutation $\sigma^* \in \Sigma_m$ such that $r_{\sigma^*(i)}\neq r_{\sigma^*(i+1)}$ for $1\leq i\leq m-1$ and $r_{\sigma^*(m)} \neq r_{\sigma^*(1)}$.
In the case where $m=2m'$ is even and $m'$ of the $r_i$ are equal, any permutation $\sigma^*$ such that $r_{\sigma^*(1)} = r_{\sigma^*(3)} = \cdots = r_{\sigma^*(2m'-1)}$ satisfies this condition. 
In all other cases let us assume that the $r_i$ are sorted by equality :
if $i<j<k$ and $r_i=r_k$ then $r_i=r_j=r_k$. If $m=2m'$ is even, we set $\sigma^*(2i-1) = i$ and $\sigma^*(2i) = m'+i$, $1\leq i\leq m'$. If $m=2m'+1$ is odd we set $\sigma^*(2i) = i$, $1\leq i\leq m'$ and $\sigma^*(2i-1) = m'+i$, $1\leq i\leq m'+1$. For example, $\sigma^* = (4\ 1\ 5\ 2\ 6\ 3\ 7)$
when $m=7$ and  $\sigma^* =(1\ 5\ 2\ 6\ 3\ 7\ 4\ 8)$ when  $m=8$.
With such a construction, we find that $|\sigma^*(i) -\sigma^*(i+1)|\geq m'$
if $m$ is odd and $|\sigma^*(i) -\sigma^*(i+1)|\geq m'-1$ if $m$ is even,
for all $1\leq i\leq m$, where we have set $\sigma^*(m+1):=\sigma^*(1)$.
Hence $\sigma\*$ satisfies the required condition, and therefore $\cyc(\lambda,r_{\sigma^*(1)},\cdots,r_{\sigma^*(m)})\neq 0$.

It is well known that if $k$ complex numbers $\alpha_1,\cdots,\alpha_k\in\C$ are such that $\alpha_1^d+\cdots+\alpha_k^d=0$, for all $1\leq d \leq k$, then $\alpha_1=\cdots =\alpha_k=0$.
Applying this property to the $m!$ complex numbers $\cyc(\lambda,r_{\sigma(1)},\cdots,r_{\sigma(m)})$, $\sigma\in\Sigma_m$, and noticing that the term corresponding to $\sigma^*$ is non zero, we see that there exists $1\leq d\leq m!$ such that $\b Q_d(\pi) = Q_d(\lambda,r_1,\cdots,r_m)\neq 0$. Since 
$\b Q_d$ has real coefficients, the numbers $\b Q_d(\pi)$ are real. Since the exponent $r_m/\deg \b Q_d$ is even it follows that $\Kpol_m(\pi)>0$, which concludes the proof of this lemma.
\sq
\nl
The following proposition, when applied to the function $\Keq=\sqrt[r_m]{\Kpol_m}$ concludes the proof of Theorem \ref{thequiv}.
\begin{prop}
\label{propequiv}
Let $m\geq 2$, and let $\Keq : \H_m\to \R_+$ be a continuous function obeying the following properties
\begin{enumerate}
\item \emph{Invariance property} : $\Keq(\pi\circ\phi) = |\det\phi|^{\frac m 2} \Keq(\pi)$. 
\item \emph{Vanishing property} : for all $\pi\in\H_m$, if $\Keq(\pi)=0$ then $K_m(\pi) = 0$.
\end{enumerate}
Then there exists a constant $C>0$ such that $\frac 1 C \Keq\leq K_m\leq C \Keq$ on $\H_m$.
\end{prop}

\proof
We first remark that $\Keq$ is homogeneous in a similar way as 
$K_m$: if $\lambda\geq 0$, then applying the invariance property to 
$\phi = \lambda^{\frac 1 m}\Id$ yields $\Keq(\pi\circ(\lambda^{\frac 1 m}\Id)) = \Keq(\lambda\pi)$ and $|\det\phi|^{\frac m 2} = \lambda$. Hence $\Keq(\lambda\pi) = \lambda \Keq(\pi)$.

Our next remark is that a converse of the vanishing property
holds: if $K_m(\pi) = 0$, then there exists a sequence $\phi_n$ of linear changes of coordinates, $\det\phi_n=1$, such that $\pi\circ\phi_n\to 0$ as $n\to\infty$. Hence $\Keq(\pi) = \Keq(\pi\circ\phi_n) \to \Keq(0)$. Furthermore, $\Keq(0) =  0$ by homogeneity. Hence $\Keq(\pi) = 0$.

We define the set $\NF_m:=\{\pi\in \H_m\sep K_m(\pi) = 0\}$. 
We also define a set $A_m\subset \H_m$ by a property ``opposite'' to the property defining $\NF_m$. A polynomial $\pi\in \H_m$ belongs to $A_m$ if and only if
$$
\|\pi\| \leq \|\pi\circ \phi\|\text{ for all } \phi \text{ such that } \det \phi = 1.
$$
The sets $\NF_m$ and $A_m$ are closed by construction, and clearly $\NF_m\cap A_m = \{0\}$. 
We now define 
$$
\underline{K_m}(\pi) = \lim_{r\to 0} \inf_{\|\pi'-\pi\|\leq r} K_m(\pi')
$$
the lower semi-continuous envelope of $K_m$. If $\underline{K_m}(\pi) = 0$ then there exists a converging sequence $\pi_n\to \pi$ such that $K_m(\pi_n)\to 0$. According to Proposition \ref{propsemicont}, it follows that $K_m(\pi)=0$ and hence $\pi\in\NF_m$. Therefore the lower semi continuous function $\underline{K_m}$ and the continuous function $\Keq$ are bounded below by a positive constant on the compact set $\{\pi\in A_m, \|\pi\|=1\}$. Since
in addition $\Keq$ is continuous and $K_m$ is upper semi-continuous, we find that
the constant
$$
C = \sup_{\pi\in A_m, \|\pi\|=1} \max\left\{ \frac{\Keq(\pi)}{\underline{K_m}(\pi)},\; \frac{K_m(\pi)}{\Keq(\pi)} \right\},
$$
is finite. By homogeneity of $K_m$ and $\Keq$, we infer that on $A_m$
\be
\label{equivAm}
\frac 1 C \Keq \leq \underline{K_m} \leq K_m \leq C \Keq.
\ee
Now, for any $\pi \in \H_m$, we consider $\hat \pi$ of minimal norm in the
closure of the set 
$\{\pi \circ \phi \sep \det \phi = 1\}$. By construction, we have $\hat \pi \in A_m$, and there exists a sequence $\phi_n$, $\det \phi_n=1$ such that $\pi \circ \phi_n \to \hat \pi$ as $n\to \infty$.
If $\hat \pi = 0$, then $K_m(\pi) = \Keq(\pi) = 0$.
Otherwise, we observe that 
$$
\underline{K_m}(\hat \pi)\leq K_m(\pi) \leq K_m(\hat \pi) \text{ and } \Keq(\hat \pi) = \Keq(\pi).
$$
Where we used the fact that $\underline{K_m}$, $K_m$ and $\Keq$ are respectively lower semi continuous, upper semi continuous, and continuous on $\H_m$.
Combining this with inequality \iref{equivAm} concludes the proof.
\sq

A natural question is to find the polynomial of smallest degree satisfying Theorem \ref{thequiv}.
This leads us to the theory of {\it invariant polynomials} 
introduced by Hilbert \cite{Hilbert} (we also refer to \cite{dix} for a survey on this subject). 
A polynomial $R$ on $\H_m$ is said to be invariant if $\mu = \frac{m\deg R} 2$ is a positive integer and for all $\pi\in \H_m$ and linear change of coordinates $\phi$, one has
\be
R(\pi\circ\phi) = (\det\phi)^\mu R(\pi).
\label{invarpolR}
\ee
We have seen for instance that $\Kpol_m$ and $\b Q_d$ are ``invariant polynomials'' on $\H_m$.

Nearly all the literature on invariant polynomials is concerned with the case of complex coefficients, both for the polynomials and the changes of variables. It is known in particular \cite{dix} that for all $m\geq 3$, there exists $m-2$ invariant polynomials $R_1,\cdots R_{m-2}$ on $\H_m$, such that for any $\pi$ (complex coefficients are allowed) and any other invariant polynomial $R$ on $\H_m$,
\be
\text{If } R_1(\pi) = \cdots = R_{m-2}(\pi) = 0,\text{ then } R(\pi) = 0.
\label{generators2}
\ee
A list of such polynomials with minimal degree is known explicitly at least when $m\leq 8$. 
Defining $r=2{\rm lcm} (\deg R_i)$ and $\Keq := \sqrt[r]{\sum_{i=1}^{m-2} \b R_i^{\frac{r}{\deg R_i} }}$, we see that $\Keq(\pi)=0$ implies $\Kpol_m(\pi) = 0$ and hence $K_m(\pi)=0$. According to proposition \ref{propequiv}, we have constructed a new, possibly simpler, equivalent of $K_m$.

For example when $m=2$ the list $(R_i)$ is reduced to the polynomial $\det$, and for $m=3$ to the polynomial $\disc$. 
For $m=4$, given $\pi = ax^4+4 b x^3 y+6c x^2 y^2+4 d x y^3+ey^4$, the list consists of the two polynomials 
$$
I = ae-4bd+3c^2,\quad J=\left|\begin{array}{ccc} a&b&c\\ b&c&d\\ c&d&e\end{array}\right|,
$$
therefore $K_4(\pi)$ is equivalent to the quantity $\sqrt[6]{|I(\pi)|^3+J(\pi)^2}$.
As $m$ increases these polynomials unfortunately become more and more complicated, and their number $m-2$ obviously increases. According to \cite{dix}, for $m=5$ the list consists of three polynomials of 
degrees $4,8,12$, while for $m=6$ it consists of $4$ polynomials of degrees $2,4,6,10$.

\section{Extension to higher dimension}
The function $K_{m,p}$  can be generalized to higher dimension $d>2$ in the following way.
We denote by $\H_{m,d}$ the set of homogeneous polynomials of degree $m$ in $d$ variables. For all $d$-dimensional simplex $T$, we define the interpolation operator $I_{m,T}$ acting from $C^0(T)$ onto the space $\P_{m-1,d}$ of polynomials of total degree $m-1$ in $d$ variables. This operator is defined by the conditions $I_{mT} v(\gamma) = v(\gamma)$ for all point $\gamma \in T$ with barycentric coordinates in the set $\{0,\frac 1 {m-1},\frac 2 {m-1},\cdots,1\}$. Following Section \S1.2, and generalizing Definition \iref{shapefunction}, we define the local interpolation error on a simplex, the global interpolation error on a mesh,
as well as the shape function.

For all $\pi \in \H_{m,d}$,
$$
K_{m,p,d}(\pi) := \inf_{|T|=1} \|\pi -I_{m,T}\pi\|_p.
$$
where the infimum is taken on all $d$-dimensional simplexes $T$ of volume $1$.
The variant $K_m^\cE$ introduced in \iref{defKE} also generalizes in higher dimension, and was introduced by Weiming Cao in \cite{C3}. Denoting by $\cE_d$ the set of $d$-dimensional ellipsoids, 
we define 
$$
K_{m,d}^\cE(\pi) = \left( \sup_{E\in \cE_d, E\subset \Lambda_\pi} |E| \right)^{-\frac {m} d},
$$
with $\Lambda_\pi = \{z\in \R^d \sep |\pi(z)|\leq 1\}$.
Similarly to Proposition \ref{propequivEllTri}, it is not hard to show that the functions $K_{m,p,d}(\pi)$ and $K_{m,d}^\cE(\pi)$ are equivalent: there exists constants $0<c\leq C$ depending only on $m,d$, such that
$$
c K_{m,d}^\cE \leq K_{m,p,d} \leq C K_{m,d}^\cE.
$$
Let $(\cT_n)_{n\geq 0}$ be a sequence of simplicial meshes (triangles if $d=2$, tetrahedrons
if $d=3$, \ldots) of a $d$-dimensional, polygonal open set $\Omega$. Generalizing 
\iref{admissibilitycond}, we say that $(\cT_n)_{n\geq 0}$ is admissible if there exists a constant $C_A$ verifying
$$
\sup_{T\in \cT_n} \diam(T) \leq C_A N^{-1/d}.
$$
The lower estimate in Theorem \ref{optitheorem} can be generalized, with straightforward adaptations in the proof. If $f\in C^m(\Omega)$ and $\seqT$ is an admissible sequence of triangulations, then
$$
\liminf_{N\to\infty} N^{\frac m d} e_{m,\cT_N}(f)_p \geq \left\|K_{m,d,p}\left(\frac{d^m f}{m!}\right)\right\|_{L^q(\Omega)}.
$$
Where $\frac 1 q:= \frac m d+\frac 1 p$.

The upper estimate in Theorem \ref{optitheorem} 
however does not generalize. The reason is that we used 
in its proof a tiling of the plane consisting of translates of a single triangle 
and of its symmetric with respect to the origin. This construction is not possible anymore in higher dimension, for example it is well known that one cannot tile the space $\R^3$, with equilateral tetrahedra. 

The generalization of the second part of Theorem \iref{optitheorem} is therefore the following.
For all $m$ and $d$, there exists a constant $C=C(m,d)>0$, such that for any polygonal open set
$\Omega \subset \R^d$ and $f\in C^m(\Omega)$ the following holds: for all $\e>0$, there exists an admissible sequence $\cT_n$ of triangulations of $\Omega$ such that  
$$
\limsup_{N\to\infty} N^{\frac m d} e_{m,\cT_N}(f)_p \leq C\left\|K_{m,d,p}\left(\frac{d^m f}{m!}\right)\right\|_{L^q(\Omega)}+\e.
$$ 
The ``tightness'' Theorem \ref{optitheorem} is partially lost due to the constant $C$. This upper bound is not new, and can be found in \cite{C3}. 
In the proof of the bidimensional theorem we define by \iref{defTiling} a tiling $\cP_R$ of the plane made of a triangle $T_R$, and some of its translates and of their symmetry with respect to the origin. In dimension $d$, the tiling $\cP_R$ cannot be constructed by the same procedure. 
The idea of the proof is to first consider a fixed tiling $\cP_0$ of the space, constituted of simplices
bounded diameter, and of volume bounded below by a positive constant, as well as a reference equilateral simplex $\Teq$ of volume $1$. We then set $\cP_R = \phi(\cP_0)$, where $\phi$ is a linear change of coordinates such that $T_R=\phi(\Teq)$. This procedure can be applied in any dimension, and yields all subsequent estimates ``up to a multiplicative constant'', which concludes the proof.

Since this upper bound is not tight anymore, and since the functions $K_{m,p,d}$ are all equivalent to $K_{m,d}^\cE$ as $p$ varies (with equivalence constants independent of $p$), there is no real need to keep track of the exponent $p$. We therefore denote by $K_{m,d}$ the function $K_{m,\infty,d}$.

For practical as well as theoretical purposes, it is desirable to have an efficient way to compute the shape function $K_{m,d}$, and an efficient algorithm to produce adapted triangulations.
The case $m=2$, which corresponds to piecewise linear elements, has been extensively studied see for instance \cite{B,CSX}. In that case there exists constants $0<c<C$, depending only on $d$, such that for all $\pi\in \H_{2,d}$,
$$
c\sqrt[d]{|\det \pi|} \leq K_{2,d}(\pi) \leq C\sqrt[d]{|\det\pi|}.
$$
where $\det \pi$ denotes the determinant of the symmetric matrix associated to $\pi$.
Furthermore, similarly to Proposition \ref{propEllipseMax}, the optimal metric for mesh refinement is given by the absolute value of the matrix of second derivatives, see \cite{B,CSX}, which is constructed in a similar way as in dimension $d=2$: with $U$ and $D=\diag(\lambda_1,\cdots,\lambda_d)$
the orthogonal and diagonal matrices such that $[\pi] = U^\trans D U$ and 
with $|D|:=\diag(|\lambda_1|,\cdots,|\lambda_d|)$, we set $h_\pi = U^\trans |D| U$. It can be shown that the matrix $h_\pi$ defines an ellipsoid of maximal volume included
into the set $\Lambda_\pi$. The case $m=2$ can therefore be regarded as solved.

For values $(m,d)$ both larger than $2$, the question of computing the shape function as well as the optimal metric is much more difficult, but we have partial answers, in particular for quadratic elements in dimension $3$. Following \S5, we need fundamental results from the theory of invariant polynomials, developed in particular by Hilbert \cite{Hilbert}. In order to apply these
results to our particular setting, we need to introduce a compatibility condition
between the degree $m$ and the dimension $d$.

\begin{definition}
We call the pair of numbers $m\geq 2$ and $d\geq 2$ ``compatible'' if and only if the following holds.
For all $\pi \in \H_{m,d}$ such that there exists a sequence $(\phi_n)_{n\geq 0}$ of $d\times d$ matrices with \emph{complex} coefficients, verifying $\det \phi_n=1$ and $\lim_{n\to\infty} \pi\circ \phi_n = 0$,
there also exists a sequence $\psi_n$ of $d\times d$ matrices with \emph{real} coefficients, verifying $\det \psi_n=1$ and $\lim_{n\to\infty} \pi\circ \psi_n = 0$.
\end{definition}

Following Hilbert \cite{Hilbert}, we say that
a polynomial $Q$ of degree $r$ defined on $\H_{m,d}$ is 
\emph{invariant} if  $\mu=\frac{mr} d$ is a positive integer and 
if for all $\pi\in\H_{m,d}$ and all linear changes of coordinates $\phi$,
\be
Q(\pi\circ\phi) = (\det \phi)^\mu Q(\pi).
\label{invpolQ}
\ee
This is a generalization of \iref{invarpolR}. We denote by $\I_{m,d}$ the set of invariant 
polynomials on $\H_{m,d}$. It is easy to see that if $\pi\in \H_{m,d}$
is such that $K_{m,d}(\pi)=0$, then $Q(\pi)=0$ for all
$Q\in \I_{m,d}$. Indeed, as seen in the proof of Proposition
\ref{vanishprop}, if $K_{m,d}(\pi)=0$ then there exists a sequence $\phi_n$ such that
$\det \phi_n=1$ and $\pi\circ\phi_n\to 0$. Therefore \iref{invpolQ}
implies that $Q(\pi)=0$. 
The following lemma shows that
the compatibility condition for the pair $(m,d)$ 
is equivalent to a converse of this property.

\begin{lemma}
\label{lemmacomp}
The pair $(m,d)$ is compatible if and only if for all $\pi \in \H_{m,d}$
$$
K_{m,d}(\pi) = 0 \text{ if and only if } Q(\pi) = 0 \text{ for all } Q\in\I_{m,d}. 
$$
\end{lemma}

\proof
We first assume that the pair $(m,d)$ is not compatible. Then there exists a polynomial $\pi_0\in\H_{m,d}$ such that there exists a sequence $\phi_n$, $\det \phi_n=1$ of matrices with \emph{complex} coefficients such that $\pi\circ\phi_n \to 0$, but there exists no such sequence with \emph{real} coefficients. This last property indicates that $K_{m,d}(\pi)>0$. On the contrary let $Q\in \I_{m,d}$ be an invariant polynomial, and set $\mu = \frac{m\deg Q} d$. The identity
$$
Q(\pi_0\circ\phi) =(\det\phi)^\mu Q(\pi_0)
$$
is valid for all $\phi$ with real coefficients, and is a \emph{polynomial identity} in the coefficients of $\phi$. Therefore it remains valid if $\phi$ has complex coefficients. If follows that $Q(\pi_0) = Q(\pi_0\circ\phi_n)$ for all $n$, and therefore $Q(\pi_0)=0$, which concludes the proof in the case where the pair $(m,d)$ is not compatible.

We now consider a compatible pair $(m,d)$.
Following Hilbert \cite{Hilbert}, we say that a polynomial $\pi\in H_{m,d}$ is a {\it null form} if and only if there exists a sequence of matrices $\phi_n$ with \emph{complex} coefficients such that $\det \phi_n=1$
and $\pi\circ \phi_n \to 0$. We denote by $\NF_{m,d}$ the set of such polynomials.
Since the pair $(m,d)$ is compatible, note that $\pi\in\NF_{m,d}$ if and only if there exists 
a sequence $\phi_n$ of matrices with \emph{real} coefficients such that $\det \phi_n=1$ and $\pi\circ\phi_n\to 0$. Hence, we find that
$$
\NF_{m,d}=\{\pi\in\H_{m,d}\sep K_{m,d}(\pi) = 0\}.
$$
Denoting by $\I_{m,d}^\C$ the set of invariant 
polynomials on $\H_{m,d}$ with {\it complex} coefficients, 
a difficult theorem of \cite{Hilbert} states that 
$$
\NF_{m,d} = \{\pi\in\H_{m,d}\sep Q(\pi) =0 \text{ for all } Q \in \I_{m,d}^\C\}
$$
It is not difficult to check that if $Q = Q_1+ i  Q_2$ where $Q_1$ and $Q_2$ have real coefficients
then \iref{invpolQ} holds for $Q$ if and only if it holds for both $Q_1$ and $Q_2$,
i.e. $Q_1$ and $Q_2$ are also invariant polynomials. 
Hence denoting by $\I_{m,d}$ the set of invariant polynomials on $\H_{m,d}$ with real coefficients, we have obtained that 
$$
\NF_{m,d} = \{\pi\in\H_{m,d}\sep Q(\pi) =0  \text{ for all }  Q \in \I_{m,d}\}
$$
which concludes the proof.
\sq

\begin{theorem}
\label{ThEquivMd}
If the pair $(m,d)$ is compatible, then there exists a polynomial $\Kpol$ on $\H_{m,d}$ (we set $r=\deg \Kpol$) and a constant $C>0$ such that for all $\pi \in \H_{m,d}$
\be 
\frac 1 C \sqrt[r]{\Kpol(\pi)}\leq K_{m,d}(\pi) \leq C\sqrt[r]{\Kpol(\pi)}.
\label{equivKmd}
\ee
If the pair $(m,d)$ is not compatible, then there does not exist such a polynomial $\Kpol$.
\end{theorem}

\proof
The proof of the non-existence property when the pair $(m,d)$ is not compatible is reported in the appendix.
Assume that the pair $(m,d)$ is compatible. We follow a reasoning very similar to \S5 to prove the equivalence \iref{equivKmd}.

We use the notations of Lemma \ref{lemmacomp} and consider the set 
$$
\NF_{m,d} = \{\pi\in\H_{m,d}\sep K_{m,d}(\pi) = 0\} = \{\pi\in\H_{m,d}\sep Q(\pi) =0, \; Q \in \I_{m,d}\}. 
$$
The ring of polynomials on a field is known to be Noetherian. 
This implies that there exists a finite family $Q_1,\cdots,Q_s\in \I_{m,d}$ of invariant polynomials on $\H_{m,d}$ such that any invariant polynomial is of the form $\sum P_i Q_i$ where 
$P_i$ are polynomials on $\H_{m,d}$. We therefore obtain
$$
\NF_{m,d} = \{\pi\in\H_{m,d}\sep Q_1(\pi) =\cdots =Q_s(\pi) = 0\}.
$$
which is a generalization of \iref{generators2}, however with no clear bound on $s$.

We now fix such a set of polynomials, set $r:= 2 {\rm lcm}_{1\leq i\leq s} \deg Q_i$, and define 
$$
\Kpol = \sum_{i=1}^s \b Q_i^{\frac r {\deg \b Q_i}}\;\;{\rm and}\;\;\Keq := \sqrt[r]{\Kpol}.
$$
Clearly $\Kpol$ is an invariant polynomial on $\H_{m,d}$, and $\NF_{m,d} = \{\pi\in\H_{m,d}\sep \Kpol(\pi)=0\}$. 
Hence the function $\Keq$ is \emph{continuous} on $\H_{m,d}$, obeys the invariance property $\Keq(\pi\circ\phi) = |\det\phi| \Keq(\pi)$, and for all $\pi\in\H_m$, $\Keq(\pi)=0$ implies $\Kpol(\pi)=0$ and therefore $K_{m,d}(\pi)=0$.
We recognize here the hypotheses of Proposition \ref{propequiv}, except that the dimension $d$ has changed.
Inspection of the proof of Proposition \ref{propequiv} shows that we use only once the fact that $d=2$, when we refer to Proposition \ref{propsemicont} and state that if $(\pi_n) \in\H_m$, $\pi_n\to \pi$ and $K_m(\pi_n)\to 0$, then $K_m(\pi)=0$. 
This property also applies to $K_{m,d}$, when the pair $(m,d)$ is compatible. Assume that $(\pi_n)\in\H_{m,d}$, $\pi_n\to \pi$ and that $K_{m,d}(\pi_n)\to 0$. Then there exists a sequence of linear changes of coordinates $\phi_n$, $\det\phi_n=1$, such that $\pi_n\circ\phi_n\to 0$. Therefore
$$
\Kpol(\pi) = \lim_{n\to \infty} \Kpol(\pi_n) = \lim_{n\to\infty} \Kpol(\pi_n\circ \phi_n) = 0
$$
It follows that $\pi\in\NF_{m,d}$, and therefore $K_{m,d}(\pi)=0$. 
Since the rest of the proof of Proposition \ref{propequiv} never uses that $d=2$, this concludes the proof of Equivalence \iref{equivKmd}.
\sq

Hence there exists a ``simple'' equivalent of $K_{m,d}$ for all compatible pairs $(m,d)$, 
while equivalents of $K_{m,d}$ for incompatible pairs need to be more sophisticated, 
or at least different from the root of a polynomial. This theorem leaves open several questions. 
The first one is to identify the list of compatible pairs $(m,d)$. It is easily shown that the pairs  $(m,2)$, $m\geq 2$, and $(2,d)$, $d\geq 2$ are compatible, but this does not provide any new results since we already derived equivalents of the shape function in these cases. More interestingly, we show in the next corollary that the pair $(3,3)$ is compatible, which corresponds to approximation by quadratic elements in dimension $3$.
There exists two generators $S$ and $T$ 
of $I_{3,3}$, which expressions are given in \cite{Salmon}
and which have respectively degree $4$ and $6$.

\begin{corol}
$\sqrt[6]{|S|^3+T^2}$ is equivalent to $K_{3,3}$ on $\H_{3,3}$.
\end{corol}

\proof
The invariants $S$ and $T$ obey the invariance properties $S(\pi\circ\phi) = (\det\phi)^4 S(\pi)$ and $T(\pi\circ\phi) = (\det\phi)^6 T(\pi)$. We intend to show that if $\pi\in\H_{3,3}$ and $S(\pi) = T(\pi) = 0$ then $K_{3,3}(\pi) = 0$. Let us first admit this property and see how to conclude the proof of this corollary. 
According to Lemma \ref{lemmacomp} the pair $(3,3)$ is compatible.
The function $\Keq := \sqrt[6]{|S|^3+T^2}$ is continuous on $\H_{3,3}$, obeys the invariance property 
$\Keq(\pi\circ\phi) = |\det\phi|\Keq(\pi)$ 
and is such that $\Keq(\pi) = 0$ implies $K_{3,3}(\pi) = 0$.
We have seen in the proof of Theorem \ref{ThEquivMd} that these properties imply the desired equivalence of $\Keq$ and $K_{3,3}$.

We now show that $S(\pi) = T(\pi) = 0$ implies $K_{3,3}(\pi) = 0$.
A polynomial $\pi\in\H_{3,3}$ can be of two types. Either it is \emph{reducible}, meaning that there exists $\pi_1\in \H_{1,3}$ (linear) and $\pi_2\in \H_{2,3}$ (quadratic) such that $\pi = \pi_1 \pi_2$, or it is \emph{irreducible}.
In the latter case according to \cite{Hartshorne}, there exists a linear change of coordinates $\phi$ and two reals $a,b$ such that 
$$
\pi\circ\phi = y^2 z - (x^3+ 3 a x z^2 + b z^3).
$$
A direct computation from the expressions given in \cite{Salmon} shows that $S(\pi\circ\phi) = a$ and $T(\pi\circ\phi) = -4b$.  If $S(\pi)=T(\pi) = 0$ then $S(\pi\circ\phi)=T(\pi\circ\phi) = 0$ and $\pi\circ \phi = y^2 z -x^3$. Therefore for all $\lambda\neq 0$, $\pi\circ \phi(\lambda x,\lambda^2 y, \lambda^{-3} z) = \lambda y^2 z -\lambda ^3 x^3$, which tends to $0$ as $\lambda\to 0$. We easily construct from this point a sequence $\phi_n$, $\det \phi_n = 1$, such that $\pi\circ\phi_n\to 0$. Therefore $K_{3,3}(\pi) = 0$.

If $\pi$ is reducible, then $\pi = \pi_1 \pi_2$ where $\pi_1$ is linear and $\pi_2$ is quadratic. Choosing a linear change of coordinates $\phi$ such that $\pi_1\circ\phi = z$ we obtain
$$
\pi\circ\phi = 3 z (a x^2+ 2 b xy + c y^2)+ z^2 ( u x + v y + w z),
$$ 
for some constants $a,b,c,u,v,w$.
Again, a direct computation from the expressions given in \cite{Salmon} shows that $S(\pi\circ\phi) = -(ac -b^2)^2$ (and $T(\pi\circ\phi) = 8 (ac -b^2)^3$). Therefore if $S(\pi) = T(\pi)= 0$ then the quadratic function $a x^2+ 2 b xy + c y^2$ of the pair of variables $(x,y)$ is degenerate. Hence there exists a linear change of coordinates $\psi$, altering only the variables $x,y$, and reals $\mu,u',v'$ such that
$$
\pi\circ\phi\circ\psi = \mu z x^2 + z^2 ( u' x + v' y + w z).
$$
It follows that $\pi\circ\phi\circ\psi(x,\lambda^{-1} y,\lambda z)$ tends to $0$ as $\lambda\to 0$. Again, this implies that $K_{3,3}(\pi)=0$, and concludes the proof of this proposition.
\sq

We could not find any example of incompatible pair $(m,d)$, which leads us to formulate the conjecture that all pairs $(m,d)$ are compatible (hence providing ``simple'' equivalents of $K_{m,d}$ in full generality). Another even more difficult problem is to derive a polynomial $\Kpol$ of minimal degree for all couples $(m,d)$ which are compatible and of interest. 

Last but not least, efficient algorithms are needed to compute metrics, from which effective triangulations are built that yield the optimal estimates. A possibility is to follow the approach proposed in \cite{C3}, i.e. solve numerically the optimization problem
$$ 
\inf \{\det H \sep  H\in S_d^+ \text{ and } \forall z \in \R^d, \<Hz,z\> \geq |\pi(z)|^{2/m}\},
$$
which amounts to minimizing a degree $d$ polynomial under an infinite set of linear constraints. When $d>2$, this minimization problem is not quadratic which makes it rather delicate. Furthermore, numerical instabilities similar to those described in Remark \ref{overfitting} can be expected to appear.

\section{Conclusion and Perspectives}

In this paper, we have introduced asymptotic estimates for the
finite element interpolation error measured in $L^p$ when the mesh is 
optimally adapted to the interpolated function.
These estimates are asymptotically sharp for functions of two variables, see Theorem \ref{optitheorem}, and precise up to a fixed multiplicative constant in higher dimension, as described in \S 6. They involve a shape function $K_{m,p}$ (or $K_{m,d,p}$ if $d>2$) which generalizes the determinant which appears in estimates for piecewise linear interpolation \cite{CSX,B,CDHM}. This function can be explicitly computed in several cases, as shows Theorem \ref{equal23}, and has equivalents of a simple form in a number of other cases, see Theorems \ref{thequiv} and \ref{ThEquivMd}.

All our results are stated and proved for sufficiently smooth functions.
One of our future objectives is to extend these results to larger classes of functions, 
and in particular to functions exhibiting discontinuities along curves. This means that
we need to give a proper meaning to the nonlinear 
quantity $K_{m,p}\left(\frac{d^m f}{m!}\right)$
for non-smooth functions.

This paper also features a constructive algorithm (similar to \cite{B}), that produces 
triangulations obeying our sharp estimates, and is described in \S3.2. However, this algorithm
becomes asymptotically effective only for a highly refined triangulation. A more practical
way to produce quasi-optimal triangulations is to adapt them to a metric, see \cite{Bois,Shew,Peyre}. 
This approach is discussed in \S 4.2. This raises the question of generating 
the appropriate metric from the (approximate) knowledge of the derivatives of the function to be interpolated. We addressed this question in the particular case of piecewise quadratic approximation in two dimensions in Theorems \ref{propEllipseMax} and \ref{family3}.

We plan to integrate this result in the PDE solver FreeFem++ in a near future. 
Note that a Mathematica source code is already available on the web \cite{sitejm}.
We also would like to derive appropriate metrics for other settings of 
degree $m$ and dimension $d$, although, as we pointed it in 
Proposition \ref{overfitting}, this might be a rather delicate matter.

We finally remark that in many applications, one seeks for error estimates in the Sobolev norms $W^{1,p}$ (or $W^{m,p}$) 
rather than in the $L^p$ norms. Finding the optimal triangulation for such norms requires a new error analysis. 
For instance, in the survey \cite{Shew2} on piecewise linear approximation,
it is observed that the metric $h_\pi = |d^2f|$ (evoked in Equation \iref{family2}) should be replaced with $h_\pi = (d^2f)^2$ for best adaptation in $H^1$ norm. In other words, the principal axes of the positive definite matrix $h_\pi$ remain the same, but its conditioning is squared. 

\appendix
\begin{center}
    {\bf APPENDIX}
  \end{center}
\section{Proof of Proposition \ref{family3}}

We consider a fixed polynomial $\pi\in\H_3$, a parameter $\alpha>0$, and look for an ellipse  $E_{\pi,\alpha}$ of maximal volume included in the set $\alpha^{-1/2} D \cap \Lambda_\pi$. Since this set is compact, a
standard argument shows that there exists at least one such ellipse. 

If $\alpha\geq \mu_\pi$, then $\alpha^{-1/2} D \subset \Lambda_\pi$ and therefore $\alpha^{-1/2} D\cap \Lambda_\pi = \alpha^{-1/2} D$. It follows that $E_{\pi,\alpha}=\alpha^{-1/2} D$, which proves part 1. 

In the following we denote by $E'_\alpha$ the ellipse defined by the matrix \iref{hBigAlpha}. Note that any ellipse containing  $D_\pi$ and included in $\Lambda_\pi$ must be tangent to $\partial\Lambda_\pi$ at the point $z_\pi$, and hence of the form $E'_\delta$ for some $\delta>0$. Clearly $E'_\delta\subset E'_\mu$ if and only if $\delta \geq \mu$. Therefore $E'_\alpha \subset \Lambda_\pi$ if and only if $\alpha\geq \alpha_\pi$. 
Let $E$ be an arbitrary ellipse, let $D_1$ the largest disc contained in $E$, and $D_2$ the smallest disc containing $E$. Then it is not hard to check that $|E| = \sqrt{|D_1| |D_2|}$.
For any $\alpha$ verifying $\alpha_\pi\leq \alpha\leq \mu_\pi$, the ellipse $E'_\alpha$ is such that $D_1 = D_\pi$, which is the largest centered disc contained in $\Lambda_\pi$, and $D_2 = \alpha^{-1/2} D$, which corresponds to the bound $2\alpha^{-1/2}$ on the diameter of $E_{\pi,\alpha}$. 
It follows that $E'_\alpha$ is an ellipse of maximal volume 
included in $\alpha^{-1/2} D\cap \Lambda_\pi$, and this concludes the proof of part 2.

Part 4 is trivial, hence we concentrate on part $3$ and assume that  $\beta_\pi\leq \alpha\leq \alpha_\pi$.

An elementary observation is that $E_{\pi,\alpha}$ must be ``blocked with respect to rotations''. Indeed assume for contradiction that $R_\theta(E_{\pi,\alpha}) \subset \Lambda_\pi$ for $\theta \in [0,\ve]$ or $[-\ve,0]$, where we denote by $R_\theta$ the rotation of angle $\theta$. Observing that the set $\cup_{\theta \in [0,\ve]} R_\theta(E_{\pi,\alpha})$ contains an ellipse of larger area than $E_{\pi,\alpha}$ and of the same diameter, we obtain a contradiction.

In the following, we say that an ellipse $E$ is quadri-tangent to $\Lambda_\pi$, when there are at least four points of tangency between $\partial E$ and $\partial\Lambda_\pi$ (a tangency point being counted twice 
if the radii of curvature of $\partial E$ and $\partial\Lambda_\pi$ coincide at this point).

The fact that $E_{\pi,\alpha}$ is ``blocked with respect to rotations'' implies that it is either quadri-tangent to $\Lambda_\pi$ or tangent to $\partial\Lambda_\pi$ at the extremities of its small axis.
In the latter case the extremities of the small axis must clearly be the points $z_\pi$ and $-z_\pi$, the closest points of $\partial \Lambda_\pi$ to the origin. It follows that $E_{\pi,\alpha}$ belongs to the family $E'_\delta$, $\delta \geq \alpha_\pi$ described above, and therefore is equal to $E'_{\alpha_\pi}$ since $\alpha\leq\alpha_\pi$. But $E'_{\alpha_\pi}$ is quadri-tangent to $\Lambda_\pi$, since otherwise we would have $E'_{\alpha_\pi-\ve}\subset \Lambda_\pi$ for some $\ve>0$.

We have now established that $E_{\pi,\alpha}$ is quadri-tangent to $\Lambda_\pi$ when $\beta_\pi\leq \alpha\leq \alpha_\pi$. This property is invariant by any linear change of coordinate: if an ellipse $E$ is quadri-tangent to $\Lambda_{\pi\circ\phi}$, then $\phi(E)$ is quadri-tangent to $\Lambda_\pi$. Furthermore if $E$ is defined by a 
symmetric positive definite matrix $H$, then $\phi(E)$ is defined by $(\phi^{-1})^\trans H\phi^{-1}$.
This remark leads us to identify the family of ellipses quadri-tangent to $\partial \Lambda_\pi$ when $\pi$ is among the four reference polynomials $x(x^2 - 3 y^2),\ x(x^2+3 y^2),\ x^2 y$ and $x^3$. In the case of $x^3$ there is no quadri-tangent ellipse and we have $\alpha_\pi=0$, therefore part 3 of the theorem is irrelevant.
In the three other cases, which respectively correspond to 
part 3 (i), (ii) and (iii), the quadri-tangent ellipses
are easily identified using the symmetries of these polynomials and the system of equations \iref{tangencyEqn}. 

The ellipses quadri-tangent to $x(x^2 + 3 y^2)$ are defined by matrices of the form $H_\lambda = \diag(\lambda,\frac{4+ \lambda^3}{3\lambda^2})$, where $0<\lambda\leq 2$. Note that $\det H_\lambda$ is decreasing on $(0,2^{\frac 1 3}]$ and increasing on $[2^{\frac 1 3},2]$. 
Given $\pi$ with $\disc\pi<0$, the optimization problem \iref{hConstrained}, therefore becomes
$$
\min_\lambda\{\det H_\lambda \sep (\phi_\pi^{-1})^\trans H_\lambda \phi_\pi^{-1} \geq \alpha \Id\}.
$$
If the constraint is met for $\lambda = 2^{1/3}$, we obtain $E_{\pi,\alpha}=E_\pi$ and therefore $\alpha\leq \beta_\pi$. Otherwise, using the monotonicity of $\lambda\mapsto \det H_\lambda$ on each side of its minimum $2^{\frac 1 3}$ we see that the matrix $H_\lambda - \alpha \phi_\pi^\trans \phi_\pi$ must be singular. Taking the determinant, we obtain an equation of degree $4$ from which $\lambda$ can be computed, and this concludes the proof of part 3 (i).

The ellipses quadri-tangent to $x(x^2 - 3 y^2)$ are defined by  $H_{\lambda, V} = V^\trans \diag(\lambda,\frac{4- \lambda^3}{3\lambda^2}) V$, where $0<\lambda\leq 1$ and $V$ is a rotation by $0,60$ or $120$ degrees. Since $\det H_{\lambda,V}$ is a decreasing function of $\lambda$ on $(0,1]$, we can apply the same reasoning as above to polynomials $\pi$ such that $\disc \pi>0$. This concludes the proof of part 3 (ii). 

Last, the ellipses quadri-tangent to $x y ^2$ are defined by $H_\lambda = \diag(\lambda,\frac 4 {27\lambda^2})$, $\lambda>0$. The determinant is a decreasing function of $\lambda$, with lower bound $0$ as $\lambda\to\infty$, and the same reasoning applies again hence concluding the proof of part 3 (iii).

\section{Proof of non existence property in Theorem \ref{ThEquivMd}}

Let $(m,d)$ be an incompatible pair. We know from Lemma \ref{lemmacomp} that there exists $\pi_0\in \H_{m,d}$ such that $K_{m,d}(\pi_0) >0$ and $Q(\pi_0)=0$ for all invariant polynomial $Q\in \I_{m,d}$.

We assume for contradiction that a polynomial $\Kpol$ satisfies inequalities \iref{equivKmd}. 
Up to replacing $\Kpol$ with $\Kpol^{2 d}$, we can assume that $\Kpol$ takes non negative values on $\H_{m,d}$ and that $\mu=\frac{mr} d$ is an integer. The rest of this proof consists in showing that $\Kpol$ \emph{needs to be} an invariant polynomial, thus leading to a contradiction since we would then
have $\Kpol(\pi_0)=0$. 
For this purpose we derive from inequalities \iref{equivKmd}, and from the invariance of $K_{m,d}$ with respect to changes of variables, the inequalities
\be
C^{-2 r} (\det \phi)^\mu\Kpol(\pi)\leq \Kpol(\pi\circ\phi)\leq C^{2 r} (\det \phi)^\mu\Kpol(\pi),
\label{equivInvarPol}
\ee
where $C$ is the constant appearing in inequalities \iref{equivKmd}.
We regard the function $\b Q(\pi,\phi) = \Kpol(\pi\circ\phi)$ as a polynomial on the vector space $V= \H_{m,d}\times M_d $, where $M_d$ denotes the space of $d\times d$ matrices, and observe that it vanishes on the \emph{hypersurface} $V_{\det} := \{(\pi,\phi)\in V \sep\det \phi = 0\}$.
Since $\phi\mapsto\det(\phi)$ is an irreducible polynomial, as shown in \cite{Bochner}, it follows that $\b Q(\pi,\phi) = (\det \phi) \b Q_1(\pi,\phi)$ for some polynomial $\b Q_1$ on $V$. Injecting this expression in inequality \iref{equivInvarPol} we obtain that $\b Q_1(\pi,\phi)$ also vanishes on the hypersurface $V_{\det}$ and the argument can be repeated. By induction we eventually obtain a polynomial $\widehat \Kpol$ on $V$ such that $\Kpol(\pi\circ\phi) = (\det \phi)^\mu \widehat \Kpol(\pi,\phi)$. It follows from inequality \iref{equivInvarPol} that for all $(\pi,\phi)\in V$
$$
C^{-2 r} \Kpol(\pi) \leq \widehat \Kpol(\pi,\phi)\leq C^{2 r} \Kpol(\pi).
$$

This implies that $\widehat \Kpol(\pi,\phi)$ does not depend on $\phi$. Otherwise, since it is a polynomial, we could find $\pi_1\in H_{m,d}$ and a sequence $\phi_n\in M_d$ such that $|\widehat \Kpol(\pi_1,\phi_n)|\to \infty$.
Therefore
$$
\Kpol(\pi\circ\phi) = (\det\phi)^\mu \widehat \Kpol(\pi,\phi)=(\det\phi)^\mu \widehat \Kpol(\pi,\Id)= (\det\phi)^\mu \Kpol(\pi).
$$
This establishes the invariance property of $\Kpol$, in contradiction with our first argument, and concludes the proof.

\subsection*{Acknowledgement}
I am grateful to to professor Nira Dyn for
her invitation in Tel Aviv University where this 
work was conceived, and to 
my PhD advisor Albert Cohen for his support
in its elaboration.

\noindent
\nl
\nl
Jean-Marie Mirebeau
\nl
UPMC Univ Paris 06, UMR 7598, Laboratoire Jacques-Louis Lions, F-75005, Paris, France
\nl
CNRS, UMR 7598, Laboratoire Jacques-Louis Lions, F-75005, Paris, France
\nl
mirebeau@ann.jussieu.fr


\begin {thebibliography} {99}

\bibitem{AF} F. Alauzet and P.J. Frey, {\it Anisotropic mesh adaptation for CFD computations}, Comput. Methods Appl. Mech. Engrg. 194, 5068-5082, 2005.
 
\bibitem{Ap} T. Apel, {\it Anisotropic finite elements: Local estimates and applications}, Advances
in Numerical Mathematics, Teubner, Stuttgart, 1999.

\bibitem{BBLS} V. Babenko, Y. Babenko, A. Ligun and A. Shumeiko,
{\it On Asymptotical
Behavior of the Optimal Linear Spline Interpolation Error of $C^2$
Functions}, East J. Approx. 12(1), 71--101, 2006.

\bibitem{B} Yuliya Babenko,  {\it Asymptotically Optimal Triangulations and Exact Asymptotics for the Optimal $L^2$-Error for Linear Spline Interpolation of $C^2$ Functions}, submitted. 

\bibitem{Bochner} M. Bocher, {\it Introduction to Higher Algebra}, Courier Dover Publications, 2004
ISBN 0486495701, 9780486495705

\bibitem{Bois} J-D. Boissonnat, C. Wormser and M. Yvinec.
{\it Locally uniform anisotropic meshing.}
To appear at the next Symposium on Computational Geometry, june 2008 (SOCG 2008) 

\bibitem{Peyre} Sebastien Bougleux and Gabriel Peyr\'e and Laurent D. Cohen. {\it Anisotropic Geodesics for Perceptual Grouping and Domain Meshing.} Proc. tenth European Conference on Computer Vision (ECCV'08), Marseille, France, October 12-18, 2008.

\bibitem{A} Y. Bourgault, M. Picasso, F. Alauzet and A. Loseille, {\it On the use of anisotropic error estimators for the adaptative solution of 3-D inviscid compressible flows,}  Int. J. Numer. Meth. Fluids. [Preprint]

\bibitem{C3} W. Cao. {\it An interpolation error estimate on anisotropic meshes in $\R^n$ and optimal metrics for mesh refinement.} SIAM J. Numer. Anal. 45 no. 6, 2368--2391, 2007.

\bibitem{C1} W. Cao, {\it Anisotropic measure of third order derivatives and the quadratic interpolation error on triangular elements}, SIAM J.Sci.Comp. 29(2007), 756-781. 

\bibitem{C2} W. Cao. {\it An interpolation error estimate in $\R^2$ based on the anisotropic measures of higher order derivatives}. Math. Comp. 77, 265-286, 2008.

\bibitem{CSX} L. Chen, P. Sun and J. Xu, {\it Optimal anisotropic meshes for
minimizing interpolation error in $L^p$-norm}, Math. of Comp. 76, 179--204, 2007.

\bibitem{CDHM} A. Cohen, N. Dyn, F. Hecht and J.-M. Mirebeau,
{\it Adaptive multiresolution analysis based on anisotropic triangulations},
preprint, Laboratoire J.-L.Lions, submitted 2008.

\bibitem{CMi} A. Cohen, J.-M. Mirebeau, {\it Greedy bisection generates optimally adapted triangulations},
preprint, Laboratoire J.-L.Lions, submitted 2008.

\bibitem{De} R. DeVore, {\it Nonlinear approximation},
Acta Numerica 51-150, 1998

\bibitem{dix} J. Dixmier, {\it Quelques aspects de la th\'eorie des invariants}, Gazette des Math\'ematiciens, vol. 43, pp. 39-64, January 1990.

\bibitem{FG} P.J. Frey and P.L. George, {\it Mesh generation. Application to finite elements}, 
Second edition. ISTE, London; John Wiley \& Sons, Inc., Hoboken, NJ, 2008.

\bibitem{Hartshorne} R. Hartshorne, {\it Algebraic Geometry}. New York: Springer-Verlag, 1999. 

\bibitem{Hilbert} D. Hilbert, {\it Theory of algebraic invariants}, Translated by R. C. Laubenbacher,
Cambridge University Press, 1993.

\bibitem{Shew}
F. Labelle and J. R. Shewchuk, {\it Anisotropic Voronoi Diagrams and Guaranteed-Quality Anisotropic Mesh Generation}, Proceedings of the Nineteenth AnnualSymposium on Computational Geometry, 191-200, 2003.

\bibitem{Lang} S. Lang, {\it Algebra}, 
Lang, Serge (2004), Algebra, Graduate Texts in Mathematics, 211 (Corrected fourth printing, revised third ed.), New York: Springer-Verlag, ISBN 978-0-387-95385-4

\bibitem{OST} P.J. Olver, G. Sapiro and A. Tannenbaum, {\it Affine invariant detection; edge maps, anisotropic diffusion
and active contours}, Acta Applicandae Mathematicae 59, 45-77, 1999.

\bibitem{Salmon} G. Salmon, {\it Higher plane curves}, third edition, 1879: http://www.archive.org/details/117724690

\bibitem{Shew2} J. R. ShewChuk, {\it What is a good linear finite element}: www.cs.berkeley.edu/\~{}jrs/papers/elemj.pdf

\bibitem{sitejm} A mathematica code for the map $\pi\mapsto h_\pi$:
www.ann.jussieu.fr/\~{}mirebeau/

\bibitem{Bamg} The 2-d anisotropic mesh generator BAMG: 
http://www.freefem.org/ff++/  (included in the FreeFem++ software)

\bibitem{Inria} A 3-d anisotropic mesh generator: http://www.math.u-bordeaux1.fr/\~{}dobj/logiciels/ mmg3d.php

\end{thebibliography}
\end{document}